\documentclass[11pt]{amsart}
\usepackage{amscd,amssymb,verbatim}
\theoremstyle{plain}
\newtheorem{Thm}{Theorem}[section]
\newtheorem{Cor}[Thm]{Corollary}
\newtheorem{Lem}[Thm]{Lemma}
\newtheorem{Prop}[Thm]{Proposition}
\newtheorem{Def}[Thm]{Definition}
\newtheorem{remark}{Remark}

\newtheorem{notation}{Notation}

\begin{document}
\title{Unconditional structures of weakly null sequences}
\author{S. A. Argyros}
\address{University of Athens}
\email{sargyros@atlas.uoa.gr}
\author{I. Gasparis}
\address{Oklahoma State University}
\email{ioagaspa@math.okstate.edu}
\date{September 06, 1998}
\begin{abstract}
The following dichotomy is established for a normalized weakly null
sequence in a Banach space: Either every subsequence admits a convex
block subsequence equivalent to the unit vector basis of \( c_0 \),
or there exists a subsequence which is boundedly convexly complete.
\end{abstract}
\maketitle
\section{Introduction} \label{S:1}
The semi-normalized weakly null sequences (i.e., sequences
\( (x_j) \) in a Banach space converging weakly to zero and
such that \( \inf_j \| x_j \| > 0 \))
being a fundamental concept in the theory of Banach spaces
have been studied extensively and several results about their
structure, and the structure of the spaces spanned by them
have been proved. We mention Bessaga-Pelczynski's theorem,
\cite{bp}, that any such sequence has a subsequence which is
Schauder basic, and the Maurey-Rosenthal's examples, 
\cite{mr}, of weakly null sequences without unconditional basic
subsequences. Both results are fundamental with enormous impact
in the development of the theory. 

After the appearance of Maurey-Rosenthal's examples, a number
of results establishing certain forms of ``restricted''
unconditionality for weakly null sequences were obtained.
We recall Elton's near unconditionality, \cite{e}, \cite{o1},
the Schreier unconditionality, stated in \cite{mr} and 
proved also later by Odell, \cite{o2}, and the Argyros-
Mercourakis - Tsarpalias
convex unconditionality, \cite{amt}.
Also, the Schreier families 
\( { \{ S_{ \xi } \} }_{ \xi < {\omega}_1 } \),
\cite{aa}, and the repeated averages hierarchy, \cite{amt}, determined
the structure of those convex combinations of a weakly null sequence
that tend to zero in norm. 

The second and the third sections of the present paper are devoted to
a unified approach of these results. Some of them are stated in a 
more general setting and the proofs, given here, are simpler than the existing
ones. 
The fourth section contains a new dichotomy for weakly null sequences.
We shall next explain our results related to this dichotomy and then 
present the results of the first two sections. We are interested in
the isomorphic structure of subsequences of a given sequence 
\( (x_j) \). Therefore, in the sequel, by a weakly null sequence
\( (x_j) \) we shall mean a normalized sequence 
which is additionally bimonotone. That is \( \| x_j \| = 1 \) and
\[                                                                                \sup_n \max \biggl \{ \biggl \| \sum_{i=1}^{n}a_i x_i \biggr \|, \quad
             \biggl \| \sum_{i=n+1}^{\infty} a_i x_i \biggr \| \biggr \}
          \leq \biggl \| \sum_{i=1}^{\infty} a_i x_i \biggr \|,
\]
for all choices of scalars \( (a_i) \).
We shall use standard Banach space facts and terminology. Throughout
this paper, \(X\) will denote a real infinite dimensional Banach space 
and \(B_X\) its closed unit ball. 
\(X^*\) stands for the the Banach space of real-valued linear
functionals on \(X\) which are continuous with respect to the norm
topology. \(c_0\) denotes the Banach space of real sequences tending to
zero, under the supremum norm. 
\(\ell_1\) is the Banach space of absolutely summable real sequences,
under the norm given by the sum of the absolute values of the
coordinates. We mention here that in the sequel, we shall often
identify the elements of \(\ell_1\) with signed measures on 
\(\mathbb{N}\).
A semi-normalized sequence \((y_j)\) in \(X\)
is called \(C\)-equivalent to the unit vector basis of
\(c_0\), if there exists a positive constant \(C\) such that
\(\| \sum_{j=1}^{n} a_j y_j \| \leq C
\max_{ j \leq n} |a_j| \),
for every \(n \in \mathbb{N}\), and all choices of scalars
\((a_j)_{j=1}^n\).
For an infinite subset \(M\) of \(\mathbb{N}\), we let
\([M]\) (resp. \([M]^{ < \infty}\)) denote 
the set of its infinite (resp. finite) subsets.

We start with some definitions and notations.

\begin{Def} \hfill
\begin{enumerate}
\item A sequence \(s=(x_j) \) in a Banach space \( X\) is said to be
{\em strongly bounded \/} if \[ \sup_n  \biggl \| 
\sum_{j=1}^{n} x_j \biggr \| \overset{ \text{def} }{=} b(s)
< \infty .\]
\item A Schauder basic sequence \(s= (x_j) \) is said to be 
{\em semi-boundedly complete }, (sb.c.), if for every sequence of coefficients
\( (a_j) \) such that \( (a_j x_j) \) is strongly bounded, we have
that \( \lim_j a_j = 0 \).
\item Let \(s= (x_j) \) be a Schauder basic sequence and 
\( \boldsymbol{a} = ( a_j ) \),
be a sequence of scalars. For every \( n \in \mathbb{N}\),
we define \[
\tau_n (\boldsymbol{a},s ) =
\sup \biggl \{ \biggl \| \sum_{j \in F }{ a_j x_j} \biggr \| \, : 
         n \leq \min F, \, \sum_{j \in F}{|a_j|} \leq 1 \biggr \},\]
and \( \tau(\boldsymbol{a},s) = \lim_n \tau_n(
\boldsymbol{a},s) \).
\end{enumerate}
\end{Def}

The basic definition related to our result is the following.

\begin{Def}
A weakly null sequence \(s= (x_j) \) is said to be 
{\em boundedly convexly complete }, (b.c.c.), if for every sequence of scalars
\(\boldsymbol{a} = (a_j) \)
such that \( (a_j x_j) \) is strongly bounded we have that
\( \tau(\boldsymbol{a},s ) = 0 \).
\end{Def} 

Evidently, a sequence \((x_j) \) is b.c.c. if and only if
the following holds 
for every scalar sequence \((a_j) \) such that
\( (a_j x_j) \) is strongly bounded: given \((F_j)\), a sequence
of consecutive finite subsets of \(\mathbb{N}\) such that
\(\inf_n{\|\sum_{j \in F_n} a_j x_j\|} > 0\), then
\(\sup_n{\sum_{j \in F_n} |a_j|} = \infty\).

It follows easily that every b.c.c. sequence is also sb.c. Also,
for a sequence \(s\) we have that
\( \tau(\boldsymbol{a},s ) > 0 \) if and only if there exists an
increasing sequence \( (F_n) \) of finite subsets of \( \mathbb{N} \)
such that 
\( \sum_{j \in F_n} |a_j| \leq 1 \),
for all \( n \in \mathbb{N} \), and
\(\inf_n \| \sum_{j \in F_n} a_j x_j \| > 0 \).

Assuming \(s\) is weakly null and that 
\(\boldsymbol{a} =(a_j) \) satisfies the stronger
condition that the series 
\( \sum_{j=1}^{\infty} a_j x_j \)
converges in norm, it does not seem clear that 
\( \tau(\boldsymbol{a},s)=0 \).
However, if \( (x_j) \) is convexly unconditional it is guaranteed that
for any such \(\boldsymbol{a} \) we have that 
\(\tau(\boldsymbol{a},s) = 0 \). The main result of the fourth
section is the following.

\begin{Thm} \label{T:1} \hfill
For every weakly null sequence \((x_j)\) one of the following two
alternatives holds exclusively:
\begin{enumerate}
\item There exists a boundedly convexly complete subsequence.
\item Every subsequence admits a convex block subsequence equivalent
to the unit vector basis of \( c_0 \).
\end{enumerate}
\end{Thm}

We recall that as a consequence of Elton's nearly unconditional
theorem, one obtains a similar dichotomy for weakly null sequences
where the two alternatives go as follows:
\begin{enumerate}
\item There exists a semi-boundedly complete subsequence.
\item Every subsequence admits a further subsequence 
equivalent to the unit vector basis
of \(c_0\).
\end{enumerate}

Thus our result may be considered as the continuation of Elton's
theorem in the direction of the deeper search in the span of the
sequence \((x_j)\), for the existence of a block subsequence equivalent
to the unit vector basis of \(c_0\). Also the alternative to the later
case is the existence of a restricted form of unconditionality
which is, in an asymptotic sense, the near unconditionality for convex
block subsequences.

Let us observe that if we assume that the sequence \((x_j)\) is
unconditional basic, then Theorem \ref{T:1} follows easily by well
known arguments. However, even in this case, our proof derives more
information on the structure of the sequence. This is a consequence
of our next result. 

\begin{Thm} \label{T:2} 
Let \(s=(x_n)\) be a weakly null sequence having no subsequence which
is b.c.c. There exist \( N \in [\mathbb{N}]\) , 
\( \xi < {\omega}_1 \)
and a constant \( C > 0 \) such that the subsequence
\((x_n)_{n \in N}\) is an
\(\ell_1^{\xi}\) spreading model, and for every \(Q \in [N] \),
\(({\xi}_{n}^{Q} \cdot s)_{n \in \mathbb{N}}\)
is \(C\)-equivalent to the unit vector basis of \(c_0\).
\end{Thm}
Where \(({\xi}_{n}^{Q} \cdot s)_{n \in \mathbb{N}} \)
is the sequence of the repeated averages of order \(\xi\)
of the sequence \((x_n)_{n \in Q}\). This concept will be explained
in the next section. 

Therefore, every weakly null sequence with no b.c.c. subsequence, 
has a subsequence which behaves
similarly to the basis of Schreier's space \(X_{\xi}\), for a countable
ordinal \(\xi\). Recall that \(X_{\xi}\) is defined as the completion
of \( ( c_{00}, \| \cdot  \|_{\xi})\) where
\[
   \|x\|_{\xi}= \sup_{F \in S_{\xi}}\sum_{n \in F}|x(n)|
\]
for \(x=(x(n)) \in c_{00} \), the space of ultimately vanishing sequences.
We do not know if the subsequence \((x_n)_{n \in M}\) resulting
from Theorem \ref{T:2} is actually equivalent to a 
subsequence of the unit vector basis
of \(X_{\xi}\) for the corresponding ordinal \(\xi\). 

The proof of Theorem \ref{T:2} which immediately implies Theorem
\ref{T:1}, is of combinatorial nature depending heavily on results
obtained in sections 2 and 3. Roughly speaking the nature of this
theorem enforces a delicate combination of the near unconditionality
with the convex unconditionality and the results related to
summability methods. 

The statement of Theorem \ref{T:1} reminds us of Rosenthal's
remarkable dichotomy, \cite{ro}, for non-trivial weak-Cauchy 
sequences. ( A weak Cauchy sequence is called non-trivial if
it is non-weakly convergent.) We recall the statement of this
theorem.
\begin{Thm} \label{hr}
Every non-trivial weak Cauchy sequence in a (real or complex)
Banach space has either a strongly summing subsequence or
a convex block basis equivalent to the summing basis.
\end{Thm}
Following \cite{ro}, a weak Cauchy basic sequence 
\((x_n)\) is said to be
{\em strongly summing} (s.s.) provided that the scalar series
\(\sum_n a_n\) converges whenever 
\(\sup_n \| \sum_{i=1}^n a_i x_i \| \)
is finite. We also recall that the basic sequence \((x_n)\) is
equivalent to the summing basis provided that for every choice
of scalars \((a_n)\), the series 
\(\sum_{n=1}^{\infty} a_n x_n \) converges if and only if the
series \(\sum_n a_n\) converges.

There are similarities but also differences between Rosenthal's
dichotomy and ours. Their relation is discussed in detail in the
last section of this paper where we also give a slight improvement of
Rosenthal's result and establish the corresponding statement to 
Elton's dichotomy for non-trivial weak Cauchy sequences.  
 
We next proceed with the results of the first two sections of this paper.
In section \ref{S:2} we present a criterion for embedding
the family \( S_{\xi}^{n}(M) \), where \( M \in [\mathbb{N}] \)
, \(\xi\) is a countable ordinal and \( n \in \mathbb{N} \), into a
hereditary family \( \mathfrak{F} \) of finite subsets of \( \mathbb{N}\).
( All relevant concepts and unexplained terminology will be thoroughly
discussed in the next section.) This criterion, Theorem 2.13, is
related to the notion of \( (\xi, M, \delta, n)\) large families
introduced in \cite{amt} and \cite{af}, and for the case \(n=1\)
it roughly says that given a hereditary family \(\mathfrak{F}\)
of finite subsets of \(\mathbb{N}\), there exists \(M \in [\mathbb{N}]\)
with \( S_{\xi}(M) \) contained in \(\mathfrak{F}\) provided that
for some subset \(\mathcal{A}\) of the probability measures on \(\mathbb{N}\)
we have that \( \sup_{F \in \mathfrak{F}} \mu(F) > \delta > 0 \),
for all \( \mu \in \mathcal{A} \), and moreover, for every
\( L \in [\mathbb{N}] \), every \( \epsilon > 0 \) and
\( \zeta < \xi \), there exists \( \mu \in \mathcal{A}\) 
supported by \(L\) and such that 
\( \sup_{F \in S_{\zeta}} \mu(F) < \epsilon \).
We apply Theorem 2.13 in section \ref{S:3} 
in order to obtain a simpler proof
(avoiding the use of the strong Cantor-Bendixson index) of
the following dichotomy established in \cite{amt}.
\begin{Thm} \label{T:3} \hfill
For a weakly null sequence \((x_n)\) in a Banach space and
\(1 \leq \xi < {\omega}_1 \), one of the following holds exclusively:
\begin{enumerate}
\item For every \( M \in [\mathbb{N}]\) there exists \( N \in [M]\)
such that for all \( L \in [N] \), \((x_n)\) is 
\( (L,\xi) \) convergent.
\item There exists \( M \in [\mathbb{N}]\), \( M= (m_i)_{i \in \mathbb{N}} \),
such that \((x_{m_i})\) is an \( l_{1}^{\xi}\) spreading model.
\end{enumerate}
\end{Thm}

It is perhaps worth noting that the summability 
methods introduced in \cite{amt} form the natural ordinal analogs
of the Cesaro summability. These methods have been already employed
in modern Banach space theory, \cite{af}, \cite{os} and it is 
possible that they can be applied to other branches of mathematics
as well.
The proof of Theorem \ref{T:3} given here is accessible to
non-specialists.  

Theorem \ref{T:3} yields
\begin{Cor} \label{C:1}
Let \((x_n)\) be a weakly null sequence in a Banach space, which is an 
\( l_{1}^{\xi}\) spreading model, for some countable
ordinal \(\xi\), yet no subsequence of \((x_n)\) is an 
\(l_{1}^{\xi + 1}\) spreading model. Then there exists a semi-normalized
convex block subsequence of \((x_n)\) which is Cesaro summable.
\end{Cor}

Let us remark here that the hypothesis of Corollary \ref{C:1} is satisfied
by any weakly null sequence in a Banach space whose Szlenk index,
\cite{s}, is equal to \({\omega}^n\), \( n \in \mathbb{N}\).
Thus we generalize the result of Alspach and Odell, \cite{ao}, who
established Corollary \ref{C:1} for weakly null sequences in
\( C({\omega}^{{\omega}^n}) \), \( n \in \mathbb{N}\). 

In the third section we give a simpler proof and a generalization of
Elton's nearly unconditional theorem, \cite{e}, \cite{o1}. More precisely
we show

\begin{Thm} \label{T:4}
Let \( s= (x_n) \) be a weakly null sequence in a Banach space and 
\(\xi\) a countable ordinal. There exists \( M \in [\mathbb{N}] \)
such that for every \(\delta \in (0,1] \)
there exists a constant \( C(\delta) > 0\)
so that the following property is satisfied: If \( L \in [M]\), 
\( n \in \mathbb{N}\) and \( (a_i)_{i=1}^{n} \) are scalars in
\([-1,1]\), then
\[\biggl \| \sum_{i \in F}a_i {\xi}_{i}^{L} \cdot s \biggr \| \leq 
C(\delta) \max \biggl \{ \delta , \, \biggl \| \sum_{i=1}^{n} 
a_i {\xi}_{i}^{L} \cdot s \biggr \| \biggr \},\] 
for all \( F \subset \{ i \leq n : \, |a_i| \geq \delta \} \).
\end{Thm}

The proof of Theorem \ref{T:4} is based on the combinatorial Lemma
3.2. The dual form of this lemma is also proved in section 3, 
Lemma 3.4, 
and as a consequence of it we obtain one of the main results of
\cite{amt} that every normalized weakly null sequence admits a
convexly unconditional subsequence. An equivalent formulation of
this result is the following:
\begin{Thm} \label{T:5}
Let \((x_n)\) be a weakly null sequence in a Banach space. There 
exists \(M \in [\mathbb{N}]\), \(M=(m_i)\),
such that for all \(\delta > 0\), there exists a constant
\(C(\delta) > 0\) so that the following property is satisfied:
If \(F \in [\mathbb{N}]^{< \infty}\) and \(({\lambda}_i)_{i \in F}\)
are scalars with \( \| \sum_{i \in F} {\lambda}_i
x_{m_i} \| > \delta \) and 
\( \sum_{i \in F} |{\lambda}_i| \leq 1 \),
then \( \| \sum_{i=1}^{\infty} a_i x_{m_i} \| > C(\delta) \),
for all choices of scalars
\((a_i)_{i=1}^{\infty} \subset c_{00} \), with \(\max_i |a_i| \leq 1\)
and such that
\(|a_i|=|{\lambda}_i|\), for all \(i \in F\).
\end{Thm} 
 
Another application of Lemma 3.4 is on the \(S_{\xi}\) unconditionality
of \(l_{1}^{\xi}\) spreading models. Recall that the sequence 
\((x_n)\) is said to be \(S_{\xi}\) unconditional, if there exists a
constant \(C > 0\) such that
\[
  \biggl \| \sum_{i \in F} a_i x_i \biggr \| \leq C
  \biggl \| \sum_{i=1}^{\infty} a_i x_i \biggr \|,
\]
for every \(F \in S_{\xi}\) and all choices of scalars
\((a_i)_{i=1}^{\infty} \subset c_{00}\).

\begin{Thm} \label{T:6}
Let \((x_n)\) be a weakly null sequence and \( \xi < {\omega}_1\).
Assume \((x_n)\) is an \(l_{1}^{\xi}\) spreading model.
There exists \(M \in [\mathbb{N}]\), \(M=(m_i)\), such that
\((x_{m_i})\) is \(S_{\xi}\) unconditional.
\end{Thm}

We also obtain, Corollary 3.7, the result on 
Schreier unconditionality, \cite{mr}, \cite{o2}, that every
normalized weakly null sequence admits, for every \(\epsilon > 0\),
a subsequence which is \(2+\epsilon\) \(S_1\) unconditional. \\
Our final results, Theorem 3.9 and Corollary 3.11, concern the
duality between \({c_0}^{\xi}\) and \(l_{1}^{\xi}\)
spreading models, and the concept of the hereditary \(\xi\)
Dunford-Pettis property. 

Before closing this section we would like to mention that according
to an unpublished result of Johnson, \cite{o1}, if every subsequence
of a normalized weakly null sequence \((x_n)\) admits a further subsequence
which is strongly bounded, then there exists a subsequence of
\((x_n)\) equivalent to the unit vector basis of \(c_0\).
Theorem \ref{T:1} immediately yields a  
generalization of Johnson's result as the following corollary
shows.
\begin{Cor} \label{C:2}
Let \((x_n)\) be a weakly null sequence in a Banach space. Assume that
every subsequence of \((x_n)\) admits a semi-normalized convex
block subsequence which is strongly bounded. Then there exists a convex
block subsequence of \((x_n)\) equivalent to the unit vector basis
of \(c_0\).
\end{Cor}
We would like to point out that even under the stronger
assumption of Corollary \ref{C:2}, the proof has to go 
through the arguments of the general case. The main difficulty
is that changing from a set \(L\) to a set \(M\), the
corresponding convex block subsequence supported by
\(L\), changes arbitrarily to a subsequence supported 
by \(M\) and thus there is no obvious way to apply
directly the infinite Ramsey's theorem.

We wish to thank H. Rosenthal for useful discussions regarding
this paper.

\section{Large families} \label{S:2}
In this section we present a criterion, Theorem 2.13, for embedding
the family \(S_{\xi}^{n}(M)\), where \(M \in [\mathbb{N}]\),
\(\xi < {\omega}_1\) and \(n \in \mathbb{N}\), 
into a hereditary family \(\mathfrak{F}\) of
finite subsets of \(\mathbb{N}\). This criterion will be applied
in section 3 in the proof of Theorem \ref{T:3}.

We shall now introduce some notation and terminology that will be
frequently used in the sequel and state all necessary definitions of
the concepts discussed in this paper. We first recall the definition
of the generalized Schreier families introduced in \cite{aa}.
It is convenient here to associate to each countable ordinal
\(\xi\), a sequence of successor ordinals \(({\xi}_n + 1)\)
such that \({\xi}_n + 1 = \xi \), for all \(n \in \mathbb{N}\),
if \(\xi\) is a successor ordinal, while \(({\xi}_n + 1)\)
strictly increases to \(\xi\), if \(\xi\) is a limit ordinal.
In the sequel we shall refer to \(({\xi}_n + 1)\) as the
sequence of ordinals associated to \(\xi\).
\begin{notation}
If \(F_1\), \(F_2\) are non-empty subsets of \(\mathbb{N}\) with
\(F_1\) finite, we denote by \(F_1 < F_2\)
the relation \(\max F_1 < \min F_2\).
If \(\mu\), \(\nu\) are finitely supported signed measures on 
\(\mathbb{N}\), we also write \(\mu < \nu \) if
\(supp \mu < supp \nu \).
\end {notation}

The Schreier families \(\{S_{\xi}\}_{\xi < {\omega}_1}\)
are defined by transfinite induction as follows:
\[  S_0 = \bigl \{ \, \{n\} : n \in \mathbb{N} \bigr \} \bigcup
    \{\emptyset\}. \]
Suppose that the families \(S_{\alpha}\)
have been defined for all \(\alpha < \xi\).

If \(\xi=\zeta + 1\), we set
\[ S_{\xi}= \biggl \{ F \in [\mathbb{N}]^{< \infty}:
   F= \bigcup_{i=1}^{n} F_i, 
   F_i \in S_{\zeta}, i \leq n,
   n \leq F_1 < \dots < F_n \biggr \} \bigcup \{\emptyset\}.\]   

If \(\xi\) is a limit ordinal, let \(({\xi}_n + 1)\)
be the sequence of ordinals associated to \(\xi\) and set
\[ S_{\xi}= \{ F \in [\mathbb{N}]^{< \infty}: n \leq \min F,
   \text{ and } F \in S_{{\xi}_n + 1}, \text{ for some }
   n \in \mathbb{N} \}.\]

\begin{Def}
A semi-normalized sequence \((x_n)\) in a Banach space is an
\(l_{1}^{\xi}\) spreading model, \(\xi < {\omega}_1\), if
there exists a constant \(C > 0\) such that
\[ \biggl \| \sum_{i \in F} a_i x_i \biggr \| \geq
   C \sum_{i \in F} |a_i| \]
for every \(F \in S_{\xi}\) and all choices of scalars
\((a_i)_{i \in F}\).
\end{Def}

We recall here that if \((x_n)\) is a sequence in some Banach space,
then the sequence \((y_n)\) is called a {\em convex block subsequence \/}
of \((x_n)\), if there exist sets \(F_i \subset \mathbb{N}\)
with \( F_1 < F_2 < \cdots \) and a sequence \((a_i)\)
of non-negative scalars such that for every \(i \in \mathbb{N}\),
\(y_i = \sum_{n \in F_i}a_n x_n \), and \(\sum_{n \in F_i}a_n =1\).
We then denote by \(suppy_i\), the support of \(y_i\), that is the
set \(\{n \in F_i: a_n > 0 \}\). We shall also adopt the notation
\(y_1 < y_2 < \cdots \) to indicate that \((y_n)\) is a convex block
subsequence of \((x_n)\).

We next pass to the definition of the repeated averages hierarchy
introduced in \cite{amt}. We let \((e_n)\) denote the unit vector
basis of \(\ell_1\). For every countable ordinal \(\xi\) and every
\(M \in [\mathbb{N}]\), we define a convex block subsequence 
\(({\xi}_{n}^{M})_{n=1}^{\infty}\)
of \((e_n)\) by transfinite induction on \(\xi\) in the following manner: \\
If \(\xi = 0\), then \({\xi}_{n}^{M} = e_{m_n}\), for all 
\(n \in \mathbb{N}\), where \(M = (m_n)\). \\
Assume that \(({\zeta}_{n}^{M})_{n=1}^{\infty}\)
has been defined for all \(\zeta < \xi\) and \(M \in [\mathbb{N}]\).
Let \(\xi = \zeta + 1 \). Set
\[ {\xi}_{1}^{M} = \frac{1}{m_1}
                   \sum_{i=1}^{m_1} {\zeta}_{i}^{M}
\]
where \(m_1 = \min M\). Suppose that
\( {\xi}_{1}^{M} < \cdots < {\xi}_{n}^{M}\)
have been defined. Let
\[
   M_n = \{ m \in M : m > \max supp {\xi}_{n}^{M} \}
\text{ and } k_n = \min M_n .
\]
Set
\[{\xi}_{n + 1}^{M} = \frac{1}{k_n} \sum_{i=1}^{k_n}
  {\zeta}_{i}^{M_n} .
\]
If \(\xi\) is a limit ordinal, let \(({\xi}_{n} + 1)\)
be the sequence of ordinals associated to \(\xi\), and let also
\(M \in [\mathbb{N}]\). Define 
\[
  {\xi}_{1}^{M} = [ {\xi}_{m_1} + 1]_{1}^{M}
\]
where \(m_1 = \min M\). Suppose that
\( {\xi}_{1}^{M} < \cdots < {\xi}_{n}^{M}\)
have been defined. Let
\[
   M_n = \{ m \in M : m > \max supp {\xi}_{n}^{M} \}
\text{ and } k_n = \min M_n .
\]
Set
\[
{\xi}_{n + 1}^{M} =  [ {\xi}_{k_n} + 1]_{1}^{M_n}.
\]     
The inductive definition of \(({\xi}_{n}^{M})_{n=1}^{\infty}\),
\(M \in [\mathbb{N}]\) is now complete.
The following properties are established in \cite{amt}. 

{\boldmath \(P1\)}: \(({\xi}_{n}^{M})_{n=1}^{\infty}\) is a convex block
subsequence of \((e_n)\) and 
\(M = \bigcup_{n=1}^{\infty} supp{\xi}_{n}^{M}\)
for all \(M \in [\mathbb{N}]\) and \(\xi < {\omega}_1\).

{\boldmath \(P2\)}: \(supp{\xi}_{n}^{M} \in S_{\xi}\), 
for all \(M \in [\mathbb{N}]\), \(\xi < {\omega}_1\) and \(n \in \mathbb{N}\).

{\boldmath \(P3\)}: If \( M, N \in [\mathbb{N}]\), \(\xi < {\omega}_1\),
and \(supp{\xi}_{i}^{M}=supp{\xi}_{i}^{N}\),
for \(i \leq k\), then \({\xi}_{i}^{M}= {\xi}_{i}^{N}\) for \(i \leq k\).

{\boldmath \(P4\)}: If \(\xi < {\omega}_1\), 
\(\{ n_k: k \in \mathbb{N}\} \subset \mathbb{N}\),
and \(\{L_k: k \in \mathbb{N}\} \subset [\mathbb{N}]\),
are such that \( supp {\xi}_{n_i}^{L_i} < supp {\xi}_{n_{i+1}}^{L_{i+1}}\),
for all \(i \in \mathbb{N}\), then letting 
\(L= \bigcup_{i=1}^{\infty} supp {\xi}_{n_i}^{L_i}\),
we have that \({\xi}_{i}^{L}={\xi}_{n_i}^{L_i}\), for all
\(i \in \mathbb{N}\).

Properties {\boldmath \(P3\)} and {\boldmath \(P4\)} are called stability
properties of the hierarchy 
 \(\{ ({\xi}_{n}^{M})_{n=1}^{\infty} : M \in [\mathbb{N}]\} \).

Next we recall the definition of \((\xi,M)\) convergent sequences.
\begin{notation}
For a sequence \(s=(x_n)\) in a Banach space and a vector 
\(\mu = \sum_{n=1}^{\infty}a_n e_n \) of \(\ell_1\), we set
\(\mu \cdot s = \sum_{n=1}^{\infty}a_n x_n\).
\end{notation}

\begin{Def}
A sequence \(s=(x_n)\) in a Banach space is called \((\xi,M)\) convergent,
\(\xi < {\omega}_1\), \(M \in [\mathbb{N}]\), if
\(\lim_{n} \|{\xi}_{n}^{M} \cdot s \|=0\).
The sequence \(s\) is called \(\xi\)-convergent, if for every 
\(M \in [\mathbb{N}]\), there exists \(N \in [M]\) such that
\(s\) is \((\xi,L)\) convergent, for all \(L \in [N]\).
\end{Def}

\begin{remark}
A sequence is \((\xi + 1 ,M)\) convergent if and only if it is
\((\xi,M)\) summable. The latter concept is defined in \cite{amt}.
Therefore the \(\xi + 1\)-convergence of a sequence is equivalent to
its \(\xi\)-summability introduced in \cite{amt}.
\end{remark}
We continue this preliminary discussion with the notion of a hereditary
family.
\begin{Def} \hfill
\begin{enumerate}
\item A family \(\mathfrak{F}\) of finite subsets of \(\mathbb{N}\)
is called {\em hereditary }, if 
for every \( G \in \mathfrak{F}\)
and \(F \subset G\) we have that \(F \in \mathfrak{F}\).
\item If \(\mathfrak{F}\) is hereditary and \(M \in [\mathbb{N}]\),
\(M=(m_i)\), then we define \(\mathfrak{F}[M]=
\{F \in \mathfrak{F}: F \subset M \}\), and
\(\mathfrak{F}(M)= \{ M(F): F \in \mathfrak{F} \}\),
where \(M(F)= \{ m_i : i \in F \} \).
\item If \(\mathfrak{F}\) is hereditary and \(n \in \mathbb{N}\),
then we set \[ {\mathfrak{F}}^n = 
\biggl \{\bigcup_{i=1}^{n}F_i: F_1 < \cdots < F_n \text{ and }
F_i \in \mathfrak{F}, i \leq n \biggr \}. \]
\end{enumerate}
\end{Def}
\begin{remark}
Evidently, if \(\mathfrak{F}\) is hereditary
then so are the families \(\mathfrak{F}[M]\),
\(\mathfrak{F}(M)\) and \({\mathfrak{F}}^n \) for every
\(n \in \mathbb{N}\) and \(M \in [\mathbb{N}]\).
It is also easily seen, by induction, that the Schreier families
\({\{S_{\xi}\}}_{\xi < {\omega}_1}\) are hereditary.
\end{remark}
We observe that for all \(\xi < {\omega}_1\) and \(M \in [\mathbb{N}]\),
\(S_{\xi}(M) \subset S_{\xi}[M]\). Note that the inverse inclusion
does not hold. However, the following result of Androulakis and Odell,
\cite{ano}, shows that the two families are in a certain sense 
comparable.

\begin{Lem} \label{T:24}
For every \(M \in [\mathbb{N}]\) there exists \(N \in [M]\) such
that for all \(\xi < {\omega}_1\),
\(\{ F \setminus \{ \min F \}: F \in S_{\xi}[N] \}
\subset S_{\xi}(M) \).
\end{Lem}

The next three lemmas describe properties of the maximal members
of \(S_{\xi}\), \(\xi < {\omega}_1\). Their proofs are easily
obtained by transfinite induction and therefore we omit them.
\begin{Lem} \label{T:25} \hfill
Let \(\xi < {\omega}_1\) and \(F \in S_{\xi}\). The following are
equivalent:
\begin{enumerate}
\item If \(F \subset G\) and \(G \in S_{\xi}\), then \(F=G\).
Thus, \(F\) is a maximal \( S_{\xi}\) set.
\item \(F \bigcup \{n\} \notin S_{\xi}\),
for all \(n \in \mathbb{N}\) with \( \max F < n\).
\item \(F \bigcup \{n\} \notin S_{\xi}\),
for some \(n \in \mathbb{N}\) with \( \max F < n\).
\end{enumerate}
\end{Lem}

\begin{Lem} \label{T:26}
Suppose that \(F_1 < \cdots < F_n \) are in \(S_{\xi}[M]\), 
\(M \in [\mathbb{N}]\). There exist \(k \leq n\) and
\(G_1 < \cdots < G_k\), maximal members of \(S_{\xi}[M]\)
with \(\min G_1 = \min F_1\) and such that
\[
\bigcup_{i=1}^{k-1}G_i \subset \bigcup_{i=1}^{n}F_i
\subset \bigcup_{i=1}^{k}G_i, \quad (G_0 = \emptyset).
\]
\end{Lem}

\begin{Lem} \label{T:27}
Let \(M \in [\mathbb{N}]\), \(\xi < {\omega}_1\).
There exists a (necesserilly) unique sequence
\(\{ F_{n}^{\xi}(M)\}_{n \in \mathbb{N}}\)
of consecutive maximal \(S_{\xi}\) sets such that
\[M= \bigcup_{n=1}^{\infty}F_{n}^{\xi}(M).\]
\end{Lem}

\begin{remark} \hfill
\begin{enumerate}
\item It is easily seen that if 
\(N=\bigcup_{n=1}^{\infty}F_{k_n}^{\xi}(M)\),
where \(k_1 < k_2 < \cdots \), then
\(F_{n}^{\xi}(N)=F_{k_n}^{\xi}(M)\),
for all \(n \in \mathbb{N}\).
\item Let \(M=(m_i)\) and \(N=(n_i)\)
be infinite subsets of \(\mathbb{N}\). Assume that for some
\(p \in \mathbb{N}\), we have that \(m_i=n_i\),
for all \(i \leq p\). If \(F_{k}^{\xi}(M)\)
is contained in \(\{m_i: i \leq p\}\),
then, \(F_{i}^{\xi}(M)=F_{i}^{\xi}(N)\),
for all \(i \leq k\).
\end{enumerate}
\end{remark}

\begin{notation}
We let \(\mathcal{M}\) denote the set of the signed measures on
\(\mathbb{N}\) whose variation does not exceed 1. Clearly,
\(\mathcal{M}\) is naturally identified with the ball of \(\ell_1\).
If \(\mu \in \mathcal{M}\) and \(\xi < {\omega}_1\), we set
\[ \|\mu\|_{\xi} = \sup \{ |\mu|(F) : F \in S_{\xi} \}. \]
\end{notation}

We would like to mention here that the proofs of the results of this 
paper rely essentially on an important principle of infinite
combinatorics known as the infinite Ramsey theorem. This theorem
was proved in several steps by Nash-Williams, \cite{nw},
Galvin and Prikry, \cite{gp} and Silver \cite{si}. Proofs
of the infinite Ramsey theorem which are not model-theoretic,
were given by Ellentuc \cite{ell}, and recently by Tsarpalias,
\cite{ts}.
We recall the statement of the theorem. \([\mathbb{N}]\) is endowed
with the topology of the pointwise convergence.
\begin{Thm} \label{T:28}
Let \(A\) be an analytic subset of \([\mathbb{N}]\). For every
\( M \in [\mathbb{N}]\) there exists \(L \in [M]\) such that
either \([L] \subset A\), or, \([L] \subset [M] \setminus A \).
\end{Thm}
In the sequel, any set satisfying the conclusion of Theorem
\ref{T:28}, will be called (completely) Ramsey.

Our next lemma is crucial for the proof of Theorem 2.13.
\begin{Lem} \label{T:29}
Let \(\mathfrak{F}\) be a hereditary family of finite subsets
of \(\mathbb{N}\), \(\xi\) a countable ordinal and \(m \in \mathbb{N}\).
For every \(M \in [\mathbb{N}]\) there exists 
\(N \in [M]\) such that either
\(\bigcup_{i=1}^{m} F_{i}^{\xi}(L) \in \mathfrak{F}\)
for all \(L \in [N]\), or,
\(\bigcup_{i=1}^{m} F_{i}^{\xi}(L) \notin \mathfrak{F}\)
for all \(L \in [N]\).
\end{Lem}
\begin{proof}
Let 
\[ \mathcal{\Delta}= \biggl \{ L \in [M]:                                                \bigcup_{i=1}^{m} F_{i}^{\xi}(L) \in \mathfrak{F} \biggr \}. \]
It follows by the second remark after lemma \ref{T:27}, that
 \(\mathcal{\Delta}\) 
is closed in \([M]\) and therefore Ramsey. Hence, there exists
\(N \in [M]\) such that either \([N] \subset \mathcal{\Delta}\),
or, \([N] \bigcap \mathcal{\Delta} = \emptyset\). 
If the former, then 
\(\bigcup_{i=1}^{m} F_{i}^{\xi}(L) \in \mathfrak{F}\),
for all \(L \in [N]\). If the latter, then
\(\bigcup_{i=1}^{m} F_{i}^{\xi}(L) \notin \mathfrak{F}\),
for all \(L \in [N]\).
\end{proof}
\begin{Def} \hfill
Suppose that \(\mathcal{A} \subset \mathcal{M}\),
\(\mathfrak{F}\) is a hereditary family of finite subsets of
\(\mathbb{N}\), \(\delta > 0\), \(M \in [\mathbb{N}]\) and 
\(\xi\) is a countable ordinal.
\begin{enumerate}
\item We shall say that \(\mathfrak{F}\) 
\(\delta\)-norms \(\mathcal{A}\), if
\( \sup_{F \in \mathfrak{F}}|\mu|(F) > \delta \),
for all \(\mu \in \mathcal{A}\).
\item \(\mathcal{A}\) is called \((\xi,M)\) large, if for 
every \(L \in [M]\), \(\zeta < \xi\) and \(\epsilon > 0\) there exists
\(\mu \in \mathcal{A}\) such that 
\(|\mu|(\mathbb{N} \setminus L) < \epsilon\)
and \(\|\mu\|_{\zeta} < \epsilon\).
\end{enumerate}
\end{Def} 

We are now ready for the proof of the main result of this section.
We first treat the case \(n=1\).
\begin{Thm} \label{T:211}
Let \(\mathfrak{F}\) be a hereditary family of finite subsets of
\(\mathbb{N}\), \(M \in [\mathbb{N}]\), \(\xi < {\omega}_1\),
and \(\delta > 0\).
Assume that there exists \(\mathcal{A} \subset \mathcal{M}\)
which is \((\xi,M)\) large and \(\delta\)-normed by
\(\mathfrak{F}\). There exists \(N \in [M]\)
such that \(S_{\xi}(N) \subset \mathfrak{F}\).
\end{Thm}
\begin{proof}
We first consider the case of \(\xi\) being a successor ordinal,
say \(\xi=\zeta +1\). Let \(P \in [M]\) and \(m \in \mathbb{N}\).
We claim that there exists \(Q \in [P]\) such that
\[ \bigcup_{i=1}^{m}F_{i}^{\zeta}(L) \in \mathfrak{F},
   \text{ for all } L \in [Q]. \]
Indeed, if this is not the case, we obtain through Lemma
\ref{T:29}, \(Q \in [P]\) such that 
\[ \bigcup_{i=1}^{m}F_{i}^{\zeta}(L) \notin \mathfrak{F},
   \text{ for all } L \in [Q]. \]
Since \(\mathcal{A}\) is \((\xi,M)\) large, there exists
\(\mu \in \mathcal{A}\) such that
\[ \|\mu\|_{\zeta} < \frac{\epsilon}{m}, \text{ and }
   |\mu|(\mathbb{N} \setminus Q ) < \frac{\epsilon}{m} \]
where \(0 < \epsilon < \frac{\delta}{2} \).
We also have that \(\mathfrak{F}\) \(\delta\) - norms \(\mathcal{A}\)
and therefore, there exists \(F_0 \in \mathfrak{F}\)
such that \(|\mu|(F_0) > \delta\).
It follows that \(|\mu|(F_0 \bigcap Q) > \delta - \epsilon \),
and thus, since \( \mathfrak{F}\) is hereditary, we can assume
that \(F_0 \subset Q\) and \(|\mu|(F_0) > \delta - \epsilon \).
Set
\[ L= F_0 \bigcup \{ q \in Q: q > \max F_0 \}\]
which belongs to \([Q]\), and choose \(k \in \mathbb{N}\)
minimal with respect to 
\(F_0 \subset \bigcup_{i=1}^{k}F_{i}^{\zeta}(L)\).
Now, \(k-1 < m\), as
\(\bigcup_{i=1}^{k-1}F_{i}^{\zeta}(L) \subset F_0\)
and so it belongs to \(\mathfrak{F}\).
Therefore \(k \leq m\) and thus
\[ \delta - \epsilon < |\mu|(F_0) \leq k \|\mu\|_{\zeta}
   < \epsilon, \text{ by the choice of } \mu, \]
which contradicts the choice of \(\epsilon\).
Hence our claim holds and we can inductively choose
\[ P_1 \supset P_2 \supset \cdots ,\text{ infinite subsets of }
M \text{ such that } \]
\[\bigcup_{i=1}^{n}F_{i}^{\zeta}(L) \in \mathfrak{F},
  \text{ for all } L \in [P_n] \text{ and } n \in \mathbb{N}.\]
Next choose \(m_1 < m_2 < \cdots\) with \(m_n \in P_n\),
for all \(n \in \mathbb{N}\).
We set \(N= (m_n)\) and claim that \(S_{\xi}(N) \subset \mathfrak{F}\).
Indeed, let \(F \in S_{\xi}\) and assume that \(\min F =n\).
Then \(F= \bigcup_{i=1}^{k}F_i\), where \(F_1 < \cdots < F_k\)
belong to \(S_{\zeta}\) and \(k \leq n\).
Applying Lemma \ref{T:26}, we obtain a finite sequence
\((G_i)_{i=1}^{l}\), \(l \leq k\), of consecutive maximal
\(S_{\zeta}\) subsets of \( N\) with \(\min G_1 = m_n\)
and such that
\[ \bigcup_{i=1}^{l-1}G_i \subset N(F) \subset \bigcup_{i=1}^{l}G_i .\]
Note that \(G_i \in [P_n]\), for \(i \leq l\) and so
there exists \(L \in [P_n]\) such that
\[ F_{i}^{\zeta}(L)= G_i, \text{ for all } i \leq l.\]
We now obtain that \(N(F) \in \mathfrak{F}\), as \(l \leq n\).
This completes the proof for the case of a successor ordinal \(\xi\).

Let now \(\xi\) be a limit ordinal and assume that the assertion of
the theorem holds for all ordinals smaller than \(\xi\).
Let \(({\xi}_n +1) \) be the sequence of ordinals associated to \(\xi\).
We can now choose by the induction hypothesis,
\[ N_1 \supset N_2 \supset \cdots \text{ infinite subsets of } M\]
such that \(S_{{\xi}_n +1}(N_n) \subset \mathfrak{F}\),
for all \(n \in \mathbb{N}\).
Suppose that \(N_i=(m_{k}^{i})_{k \in \mathbb{N}}\),
for all \( i \in \mathbb{N}\), and choose 
\(m_1 < m_2 < \cdots \) such that \(m_i \in N_i\)
and \(m_i > m_{i}^{i}\),
for all \( i \in \mathbb{N}\). Set \(N=(m_i)\) and it is easy to
see that \(S_{\xi}(N) \subset \mathfrak{F}\).
\end{proof}
\begin{notation}
Let \((x_n)\) be a sequence in a Banach space \(X\) and \(\epsilon > 0\).
We set
\[ {\mathfrak{F}}_{\epsilon}= \{ F \in [\mathbb{N}]^{< \infty}: \exists
   \, x^{*} \in B_{X^{*}} \text{ with } |x^*(x_n)| \geq \epsilon ,
   \forall \, n \in F \}. \]
Clearly, \({\mathfrak{F}}_{\epsilon}\) is a hereditary family.
\end{notation}

\begin{Cor} \label{T:212} \hfill
Let \((x_n)\) be a sequence in a Banach space \(X\) and \(\delta > 0\).
Let \(\xi\) be a countable ordinal and suppose that there exists a set
\(\mathcal{A}\) of absolutely sub-convex combinations of \((x_n)\) such that
\begin{enumerate}
\item \(\|x\| > \delta\), for all \(x \in \mathcal{A}\).
\item For every \(L \in [\mathbb{N}]\), \(\epsilon > 0\) and
\(\zeta < \xi\), there exists \(x \in \mathcal{A}\),
\(x= \sum_{i=1}^{\infty}a_i x_i\), such that
\(\sum_{i \notin L}|a_i| < \epsilon\) and
\(\sum_{i \in F}|a_i| < \epsilon\), for all \(F \in S_{\zeta}\).
\end{enumerate}
Then there exists \(N \in [\mathbb{N}]\) such that
\(S_{\xi}(N) \subset {\mathfrak{F}}_{\frac{\delta}{2}}\).
\end{Cor}

\begin{Thm} \label{T:213}
Let \(\mathfrak{F}\) be a hereditary family, \(\xi < {\omega}_1\),
\(M \in [\mathbb{N}]\),
\(\delta > 0\) and \(n \in \mathbb{N}\). Assume that there exists
\(\mathcal{A} \subset \mathcal{M}\), consisting of finitely supported
measures, which is \((\xi,M)\) large.
Assume further that if \(\mu_1 < \cdots < \mu_n\) and
\(\mu_i \in \mathcal{A}\), for \(i \leq n\), then there exists
\(G \in \mathfrak{F}\) such that
\(|\mu_i|(G) > \delta\), for all \(i \leq n\).
Then, there exists \(N \in [M]\) such that 
\(S_{\xi}^{n}(M) \subset \mathfrak{F}\).
\end{Thm}

We shall need the following
\begin{Lem} \label{T:214}
Let \(M \in [\mathbb{N}]\), \(n \in \mathbb{N}\) and
\(\xi < {\omega}_1\). Suppose that
\[ 
S_{\xi}^{n}(M) = \bigcup_{i=1}^{k} {\mathfrak{F}}_i. \]
where \({\mathfrak{F}}_i\) is hereditary for all \(i \leq k\).
Then there exist \(N \in [M]\) and \(i_0 \leq k\)
such that \(S_{\xi}^{n}(N) \subset {\mathfrak{F}}_{i_0}\).
\end{Lem}
\begin{proof}
Define
\[ 
  {\mathcal{\Delta}}_i = \biggl \{ L \in [\mathbb{N}]:
M \biggl (\bigcup_{j=1}^{n}F_{j}^{\xi}(L) \biggr ) 
\in {\mathfrak{F}}_i \biggr \},
\text{ for } i \leq k. \] 
Evidently, \({\mathcal{\Delta}}_i \) is closed in
\([\mathbb{N}]\) and therefore Ramsey, for all \(i \leq k\).
Of course, by our assumption, we have that
\[ [\mathbb{N}]= \bigcup_{i=1}^{k}{\mathcal{\Delta}}_i . \]
We may now choose \(P \in [\mathbb{N}]\) and \(i_0 \leq k\)
such that \([P] \subset {\mathcal{\Delta}}_{i_0}\).
Let \(N=M(P)\), and clearly 
\(S_{\xi}^{n}(N) \subset {\mathfrak{F}}_{i_0}\).
\end{proof}

\begin{proof}[Proof of Theorem \ref{T:213}] \hfill
By induction on \(n\). The case \(n=1\) has been settled in Theorem
\ref{T:211}. So assume \(n \geq 2\).
Choose \(\mathcal{B}\) a countable dense subset of \(\mathcal{A}\).
Clearly, \(\mathcal{B}\) satisfies the same assumptions as 
\(\mathcal{A}\) does in the hypothesis of Theorem \ref{T:213}.
Therefore, without loss of generality, we shall assume that
\(\mathcal{A}\) itself is countable. Let
\((a_k)\) be an enumeration of the elements of \(\mathcal{A}\).
We claim that there exist:
\begin{enumerate}
\item A sequence \((G_k)\) of elements of \(\mathfrak{F}\) such that
\(G_k \subset suppa_k\) and \(|a_k|(G_k) > \delta\), for all
\(k \in \mathbb{N}\).
\item A decreasing sequence \((M_k)\) of infinite subsets of \(M\)
such that \(suppa_k < M_k\) and 
\(G_k \bigcup F \in \mathfrak{F}\), for all
\(F \in S_{\xi}^{n-1}(M_k)\), and all \(k \in \mathbb{N}\).
\end{enumerate}

Indeed, suppose that \(G_1, \cdots , G_{k-1}\)
and \(M_1 \supset \cdots \supset M_{k-1}\) satisfying 1 and 2
have been constructed. Let \(G \in \mathfrak{F}\) such that
\(G \subset suppa_k\) and \(|a_k|(G) > \delta\). Define
\[ {\mathfrak{F}}_G = \{ F \in \mathfrak{F}: G < F \text{ and }
                         G \bigcup F \in  \mathfrak{F} \}. \]
Let \(\mathcal{\Delta}= \bigcup_{G} {\mathfrak{F}}_G \),
where the union is taken over all possible subsets \(G\) of
\(suppa_k\) which belong to \(\mathfrak{F}\) and satisfy the relation
\(|a_k|(G) > \delta\). 
Of course, \(\mathcal{\Delta}\) is hereditary and the hypothesis of
Theorem \ref{T:213} is satisfied by the family \(\mathcal{\Delta}\),
the integer \(n-1\) and the set of measures 
\(\{ a_i: a_k < a_i \}\).
By the induction hypothesis there exists, for all \(k \in \mathbb{N}\),
\(N_k \in [M_{k-1}]\) such that \(suppa_k < N_k\) and 
\(S_{\xi}^{n-1}(N_k) \subset \mathcal{\Delta}\).
Next choose according to Lemma \ref{T:214},
\(G_k  \in \mathfrak{F}\) and \(M_k \in [N_k]\) such that
\(G_k  \subset supp a_k\), \(|a_k|(G_k) > \delta\) and
\(S_{\xi}^{n-1}(M_k) \subset {\mathfrak{F}}_{G_k} \).
It follows now that \(G_k \bigcup F \in \mathfrak{F} \),
for all \(F \in S_{\xi}^{n-1}(M_k)\).
This completes the inductive construction and our claim holds.

Let now \(F\) be a finite subset of \(\mathbb{N}\) and set
\[ k_F = \min \{ k \in \mathbb{N}: F \subset G_k \}. \]
We let \(k_F=0\), if \(F\) is not contained in \(G_k\), for
all \(k \in \mathbb{N}\). Inductively we construct
a sequence of positive integers,
\(m_1 < m_2 < \cdots \), in the following manner: Suppose that
\(m_1< \cdots < m_d\) and \(q_1 \leq \cdots \leq q_{d-1}\)
have been constructed. ( \(m_1\) is chosen arbitrarily in
\(M_1\) and \(q_0 = 0\).) We set
\[ q_d = \max \bigl \{ k_F : F \subset \{m_1, \cdots , m_d \} \bigr \} 
\vee q_{d-1}. \]
Next choose \(m_{d+1} \in M_{q_d +1}\) such that \(m_{d+1} > m_d\)
and \(m_{d+1} > m_{d+1}^{q_d + 1}\), where
\(M_i = (m_j^i)_{j=1}^{\infty}\), for all \( i \in \mathbb{N}\).
Let \(P=(m_i)\). We now claim that if \( F \subset P\) is contained
in \(G_k\), for some \(k \in \mathbb{N}\), then 
\(F \bigcup H \in \mathfrak{F}\), for all
\(H \in S_{\xi}^{n-1}(P)\) such that \(F < H\).
Indeed, let \(m_d = \max F\). Then \(k_F \leq q_d \) and
\(F \subset G_{k_F} \). It suffices to show that 
\(H \in S_{\xi}^{n-1}(M_{k_F})\) for all 
\(H \in S_{\xi}^{n-1}(P)\) such that \(F < H\). Our claim will then follow
by condition 2 above. To this end, let
\(R \in S_{\xi}(P)\), such that \(F < R \). Suppose that
\(R=\{ m_{i_1}, \cdots , m_{i_t} \} \), where
\( \{i_1, \cdots , i_t\} \) belongs to \(S_{\xi}\).
Now, for all \(r \leq t\),
\[ m_{i_r} \geq m_{d+1} \text{ and thus } i_r \geq d+1.\]
So,
\[ m_{i_r} > m_{i_r}^{q_{i_{r} - 1}} \geq m_{i_r}^{q_d} \geq 
   m_{i_r}^{k_F}, \text{ as } k_F \leq q_d . \]
Hence, \(R \in S_{\xi}(M_{k_F})\), as claimed.
Finally, consider the hereditary family
\[ \mathfrak{D} = \bigcup_{k=1}^{\infty} 
   \bigl \{ F : F \subset G_k \bigr \} . \]
Clearly, the hypothesis of Theorem \ref{T:211} is satisfied by 
\(\mathfrak{D} \) and the set of measures
\(\{a_k | G_k : k \in \mathbb{N} \} \). 
We can thus find \(N \in [P]\) such that
\(S_{\xi}(N) \subset \mathfrak{D} \).
It is now easily verified that \( S_{\xi}^{n}(N) \subset \mathfrak{F} \).
\end{proof}

Our next proposition will enable us verify that the set of measures
\[ {\mathcal{A}}_{\xi} = \{ {\xi}_{n}^{L}: n \in \mathbb{N} ,
   L \in [\mathbb{N}] \} \] 
is \((\xi,\mathbb{N})\) large, for all \( \xi < {\omega}_1 \).
\begin{Prop} \label{T:215}
For every \( M \in [\mathbb{N}] \), \(\epsilon > 0\), and all ordinals
\( \zeta < \xi < {\omega}_1 \), there exists \( N \in [M]\)
(depending on \(M\), \(\epsilon\), \(\zeta\), \(\xi\)) 
such that \(\| {\xi}_{n}^{L} \|_{\zeta} < \epsilon\),
for all \( L \in [N] \) and \( n \in \mathbb{N} \).
\end{Prop}
\begin{proof}
It suffices to show the following:

{\em Claim: \/} Let \(\xi < {\omega}_1\). For all \(\zeta < \xi\), 
\(\epsilon > 0\) and \( M \in [\mathbb{N}] \),
there exists \(L \in [M]\) such that  
\(\| {\xi}_{1}^{L} \|_{\zeta} < \epsilon\).

Indeed, assuming our claim holds, we observe that the set
\(\{ L \in [M] : \| {\xi}_{1}^{L} \|_{\zeta} < \epsilon \} \),
is closed in \([M]\) and therefore Ramsey. Our claim now
yields the existence of \( N \in [M] \) such that
\(\| {\xi}_{1}^{L} \|_{\zeta} < \epsilon\), for every
\(L \in [N]\). By stability, we obtain the assertion of the proposition.

We shall prove our claim by transfinite induction on \(\xi\).
If \(\xi =1\), then \(\zeta=0\) and the claim is easily verified.
Assuming our claim holds for all ordinals smaller than \(\xi\),
let first \(\xi\) be a limit ordinal. Let also
\(({\xi}_n + 1)\) be the sequence of ordinals associated to \(\xi\).
Suppose now that \(\zeta < \xi\) and choose \(m \in M\) so that
\[ \zeta < {\xi}_m  \text{ and } \frac{1}{m} < \frac{\epsilon}{2}. \]
We apply the induction hypothesis on the ordinal \({\xi}_m \)
and the set \(L_1 = \{ i \in M : i \geq m \} \)
to obtain \(L_2, \cdots , L_m\),
infinite subsets of \(L_1 \) such that
\[ [{\xi}_m]_{1}^{L_1} < [{\xi}_m]_{1}^{L_2} < \cdots < 
   [{\xi}_m]_{1}^{L_m} \text{ and } \]
\[ \|[{\xi}_m]_{1}^{L_i}\|_{\zeta} < \frac{\epsilon}{2m},
   \text{ for } 2 \leq i \leq m . \]
By stability property \(P_4\), there exists \(N \in [M]\) with
\(\min N = m\) and such that
\[[{\xi}_m]_{i}^{N} = [{\xi}_m]_{1}^{L_i}, \text{ for }
  i \leq m . \]
Now, 
\[ [{\xi}_m +1 ]_{1}^{N} = \frac{1}{m} \sum_{i=1}^{m}
   [{\xi}_m]_{i}^{N} = {\xi}_{1}^{N} . \]
Hence,
\[ \| {\xi}_{1}^{N}\|_{\zeta} \leq \frac{1}{m} + \frac{m-1}{2m} 
   \epsilon < \epsilon . \]
So our claim holds if \(\xi\) is a limit ordinal. 

Suppose now that \(\xi = \alpha + 1\). If \( \zeta < \alpha\),
choose according to the induction hypothesis \(L_1 \in [M]\)
such that \( \|{\alpha}_{1}^{L_1}\|_{\zeta} < \epsilon\).
Let \(m= \min supp {\alpha}_{1}^{L_1}\)
and choose again by the induction hypothesis, 
\(L_2, \cdots , L_m \), infinite subsets of \(M\) such that
\[ {\alpha}_{1}^{L_1} < {\alpha}_{1}^{L_2} < \cdots < 
   {\alpha}_{1}^{L_m} \text{ and } \]
\[ \| {\alpha}_{1}^{L_i}\|_{\zeta} < \epsilon ,
   \text{ for } 2 \leq i \leq m . \]
But once again, by stability, there exists \(N \in [M]\) with
\(\min N = m\) and such that 
\[ {\alpha}_{i}^{N}   =  {\alpha}_{1}^{L_i},
   \text{ for } i \leq m . \]
Now, \({\xi}_{1}^{N}  = \frac{1}{m} \sum_{i=1}^{m}
       {\alpha}_{i}^{N} \),
and thus, \( \| {\xi}_{1}^{N} \|_{\zeta} < \epsilon \).

The final case to consider is when \(\zeta = \alpha\).
Let \(({\beta}_i +1)\) be the sequence of ordinals associated to
\(\alpha\). Choose \( m \in M\) such that
\(\frac{1}{m} < \frac{\epsilon}{2}\).
Set \( L_1 = \{ n \in M : n \geq m \} \).
It follows that 
\[ [{\beta}_{m} + 1 ]_{1}^{L_1} = {\alpha}_{1}^{L_1} .\]
Let \( k_1 = \max supp {\alpha}_{1}^{L_1} \).
Choose according to the induction hypothesis \(L_2 \in [M]\)
with \(k_1 < \min L_2 \) and such that
\[ \| {\alpha}_{1}^{L_2} \|_{{\beta}_j} < 
   \frac{\epsilon}{2k_1} , \text{ for all } j \leq k_1 . \]
Let \( k_2  = \max supp {\alpha}_{1}^{L_2} \).
Successive repetitions of the previous argument yield
\[ {\alpha}_{1}^{L_1} < {\alpha}_{1}^{L_2} < \cdots < 
   {\alpha}_{1}^{L_m} \text{ with }
L_i \in [M] \text{ for } i \leq m , \]
such that if \( k_i = \max supp {\alpha}_{1}^{L_i} \)
for \( i \leq m\), then
\[ 
   \| {\alpha}_{1}^{L_i} \|_{{\beta}_j} <    
  \frac{\epsilon}{2k_{i-1}} , \text{ for all } j \leq k_{i-1}
  \text{ and } 2 \leq i \leq m . \]
Stability now guarantees the existence of \(L \in [M]\)
with \(\min L =m \) and such that
\[ 
{\alpha}_{i}^{L} = {\alpha}_{1}^{L_i} , \text{ for all }
 i \leq m .\]
Now, \( {\xi}_{1}^{L} = \frac{1}{m} \sum_{i=1}^{m} {\alpha}_{i}^{L} \)
and it remains to show that \( \| {\xi}_{1}^{L} \|_{\alpha} 
< \epsilon \). 
Indeed, let \( G \in S_{\alpha}\) and choose \(i_0 \leq m \)
minimal with respect to \( G \bigcap supp {\alpha}_{i_0}^{L} 
\ne \emptyset \).
Let \(l = \min G \) and observe that \(l \leq k_{i_0}\).
Choose \(p \leq l\) such that 
\(G \in S_{{\beta}_p + 1} \).
There exist \(q \leq l\) and \((G_j)_{j=1}^{q}\)
consecutive members of \(S_{{\beta}_p}\) such that
\(G= \bigcup_{j=1}^{q} G_j\). Note that
\[ q \leq l \leq k_{i_0} \text{ and also } p \leq l \leq k_{i_0} . \]
Thus,
\[ \|{\alpha}_{i}^{L}\|_{{\beta}_p} < \frac{\epsilon}{2k_{i_0}} ,
   \text{ for } i_0 < i \leq m .\]
Therefore, 
\[ {\alpha}_{i}^{L} (G) \leq q \|{\alpha}_{i}^{L}\|_{{\beta}_p} 
< k_{i_0} \frac{\epsilon}{2k_{i_0}} = \frac{\epsilon}{2},
\text{ for } i_0 < i \leq m .\]
Hence,
\[ {\xi}_{1}^{L}(G)  = \frac{1}{m} \biggl [ {\alpha}_{i_0}^{L}(G) +
   \sum_{i=i_0 + 1}^{m} {\alpha}_{i}^{L}(G) \biggr ]
< \frac{1}{m} + \frac{\epsilon}{2} < \epsilon . \]
and so \( \| {\xi}_{1}^{L} \|_{\alpha} < \epsilon \). 
The proof of the claim is now complete.
\end{proof}
We next recall the concept of the
\((n, \xi, M , \delta)\) large families, \cite{amt}, \cite{af}.

\begin{Def}
Let \(M \in [\mathbb{N}]\), \(\xi < {\omega}_1\), \(\delta > 0\)
and \(n \in \mathbb{N}\).
The hereditary family \(\mathfrak{F}\) is called \((n, \xi, M , \delta)\)
large provided that for all \(N \in [M]\) there exists
\(F \in \mathfrak{F}\) such that
\( {\xi}_{i}^{N}(F) > \delta\), for all \(i \leq n\).
\end{Def}
\begin{Cor} \label{T:217}
Let \(\mathfrak{F}\) be \((n, \xi, M , \delta)\) large. There exists
\(N \in [M]\) such that \(S_{\xi}^{n}(N) \subset \mathfrak{F}\).
\end{Cor}
\begin{proof}
Proposition \ref{T:215} and the fact that \(\mathfrak{F}\)
is \((n, \xi, M , \delta)\) large, immediately yield that
\(\mathfrak{F}\) and \({\mathcal{A}}_{\xi}\)
satisfy the hypothesis of Theorem \ref{T:213}.
The assertion of the corollary now follows.
\end{proof}

\section{A generalization of Elton's theorem.} \label{S:3}
In the first part of this section we give the proof of Theorem \ref{T:4}
which extends Elton's nearly unconditional theorem \cite{e}, \cite{o1}.
The second part of section \ref{S:3} is devoted to the proof of the
results concerning the summability methods, the convex unconditionality
and the \(S_{\xi}\) unconditionality. We also discuss the duality
between \(c_{0}^{\xi}\) and \(l_{1}^{\xi}\) spreading models as well
as the \(\xi\) Dunford - Pettis property.

Let now \(s=(x_n)\) be a weakly null sequence in the Banach space \(X\).
Recall here that \(s\) is additionally assumed to be a normalized,
bimonotone, Schauder basic sequence. For this fixed sequence \(s\)
and the countable ordinal \(\xi\), we have the following
\begin{Def}
Let \(M \in [\mathbb{N}]\), \(\lambda > 0\) and \(k \in \mathbb{N}\).
The functional \(x^* \in B_{X^*}\) is said to be \(\xi\)-{\em good \/}
for \((\lambda, M, k)\), if 
\(x^*( {\xi}_{i}^{M} \cdot s) \geq 0 \), for all \(i \leq k\),
and \(\sum_{i=1}^{k}x^*( {\xi}_{i}^{M} \cdot s) > \lambda \).
\end{Def} 

The main tool for proving Theorem \ref{T:4} 
is the following combinatorial result.
\begin{Lem} \label{T:32}
Let \(\epsilon > 0\), \(\lambda > 0\) and \(N \in [\mathbb{N}]\).
There exists \(M \in [N]\) satisfying the following property:
If \(L \in [M]\),
\(k \in \mathbb{N}\) and there exists \(x^* \in B_{X^*}\) 
which is \(\xi\)-good for
\((\lambda , L , k)\), then there exists  \(y^* \in B_{X^*}\) 
which is \(\xi\)-good for
\((\lambda , L , k)\) and such that
\[ \sum_{n \in M \setminus \bigcup_{i=1}^{k}supp {\xi}_{i}^{L}}
   |y^*(x_n)| < \epsilon. \]
\end{Lem}
\begin{proof}
Choose first \(({\epsilon}_i)_{i=0}^{\infty}\), 
a sequence of positive scalars
such that \(\sum_{i=0}^{\infty} {\epsilon}_i < \epsilon\).
Let \(T_0 = \emptyset\) and \(T_n = \{1, \cdots , n\}\),
if \(n \in \mathbb{N}\).
By an \(n\)-tuple of positive integers 
\((m_i)_{i \in T_n}\), we shall either mean the empty
tuple, if \(n=0\), or, that \(m_1 < \cdots < m_n \),
if \(n \geq 1\).
Let now \(n \in \mathbb{N} \bigcup \{0\}\) and 
\(F \subset T_n\). 
The \(n\)-tuple \((m_i)_{i \in T_n}\)
and the infinite subset \(L\) of \(\mathbb{N}\),
(\(L=(l_i)\)), are said to
satisfy property \((F-E_n)\), provided that
\(m_n < l_1\), if \(n \geq 1\), and the following statement holds:

If \(k \in \mathbb{N}\),
and there exists \(x^* \in B_{X^*}\) which is \(\xi\)-good for
\((\lambda , \{m_j : j \in T_n \setminus F\} \bigcup
\{l_j : j \geq 2 \} , k)\), 
then there exists
\(y^* \in B_{X^*}\) which is \(\xi\)-good for
\((\lambda , \{m_j : j \in T_n \setminus F\}  \bigcup
\{l_j : j \geq 2 \} , k)\)
and such that
\[ \sum_{j \in F} |y^*(x_{m_j})| + |y^*(x_{l_1})|
< \sum_{i=0}^{n}{\epsilon}_i .\] 
Let us also say that \((m_i)_{i \in T_n}\) and \(L\)
satisfy property \((E_n)\), if they satisfy property
\((F-E_n)\), for every \(F \subset T_n \). 

We shall inductively construct an increasing sequence
\((m_n)_{n=1}^{\infty}\) 
of elements of \(N\), and a decreasing sequence
\((M_n)_{n=0}^{\infty}\) of infinite subsets of \(N\) 
with \(m_n \in M_{n-1}\), if \(n \geq 1\),
and so that 
for every \(n \in \mathbb{N}\bigcup \{0\}\), 
if \(L \in [M_n]\), then
\((m_i)_{i \in T_n}\) and \(L\) satisfy property \((E_n)\).

The first inductive step is similar to the general one and
therefore we shall not explicit it. Assume that 
\((m_i)_{i \in T_n}\) and 
\(M_0 \supset \cdots \supset M_n \), infinite subsets of \(N\),
with \(m_i \in M_{i-1}\) for \(1 \leq i \leq n\)
have been constructed so that if \(i \leq n\) and 
\(L \in [M_i]\), then  \((m_j)_{j \in T_i}\) and \(L\) 
satisfy property \((E_i)\). 
Let \(m_{n+1}= \min M_n\). Fix \(F \subset T_{n+1}\)
and define
\[ {\Delta}_F = \{ L \in [M_n], L=(l_i),: (m_i)_{i \in T_{n+1}}
   \text{ and } L \text{ satisfy } (F-E_{n+1}) \}. \]
Clearly, \({\Delta}_F\) is closed in \([M_n]\) and therefore
Ramsey. Suppose that for some \(P \in [M_n]\), \(P=(p_i)\),
we had that
\([P] \bigcap {\Delta}_F = \emptyset \).
Let \(q \in \mathbb{N}\) and set
\[ L_j = \{ p_j \} \bigcup \{ p_i : i > q \}, \text{ for all }
   j \leq q. \]
Since \(L_j \notin {\Delta}_F \), for all \(j \leq q \), 
there exist integers \((k_j)_{j=1}^{q}\)
and functionals \((x_{j}^{*})_{j=1}^{q}\) in \(B_{X^*}\)
so that letting \(R=\{m_i : i \in T_{n+1} \setminus F \} \bigcup
\{p_i : i > q \}\), we have that for all \(j \leq q \)
\[ x_{j}^{*} \text{ is } \xi \text{-good for } 
(\lambda , R , k_j ),\]
and moreover, if \(y^* \in B_{X^*}\) 
is \(\xi\)-good for \((\lambda , R , k_j )\) then,
\[  \sum_{i \in F} |y^*(x_{m_i})| + |y^*(x_{p_j})|
 \geq \sum_{i=0}^{n+1}{\epsilon}_i .\]
Next choose \(j_0 \leq q \) such that \(k_{j_0}= \min
\{k_j : j \leq q\} \).
We observe that if \(y^*\) is \(\xi\)-good for 
\((\lambda , R , k_{j_0} )\) 
and \(y^*(x_i)=0\), for all \(i > \max supp {\xi}_{k_{j_0}}^{R}\),
then \(y^*\) is also \(\xi\)-good for 
\((\lambda , R , k_j )\), for
every \( j \leq q\).

Now let \(t= \max F\) and note that
\[ R= \{m_i : i \in T_{t-1} \setminus F\} \bigcup 
   \{ m_i : t < i \leq n + 1\} \bigcup
   \{ p_i : i > q \}. \]
By the induction hypothesis, 
since \((m_i)_{i \in T_{t-1}}\) and 
\( \{ m_i : t \leq i \leq n + 1\} \bigcup
   \{ p_i : i > q \}\)
satisfy property \((E_{t-1})\),
 there exists 
\(z^* \in B_{X^*}\) \(\xi\)-good for 
\((\lambda , R , k_{j_0} )\)
and such that 
\[  \sum_{i \in F \setminus \{t\}} |z^*(x_{m_i})| + |z^*(x_{m_t})|
  < \sum_{i=0}^{t-1}{\epsilon}_i .\]
Thus, \[\sum_{i \in F} |z^*(x_{m_i})|< \sum_{i=0}^{n}{\epsilon}_i . \]
Without loss of generality, since \((x_i)\) is bimonotone,
we can assume that \(z^*(x_i)=0 \), for all
\(i > \max supp {\xi}_{k_{j_0}}^{R}\).
Our previous observation yields that 
\(z^*\) is \(\xi\)-good for 
\((\lambda , R , k_j )\), for
all \( j \leq q\), and thus,
\[  \sum_{i \in F} |z^*(x_{m_i})| + |z^*(x_{p_j})|
 \geq \sum_{i=0}^{n+1}{\epsilon}_i , \text{ for all } j \leq q .\]
Hence, \(|z^*(x_{p_j})| \geq {\epsilon}_{n+1}\), for all
\(j \leq q\).
We have reached a contradiction since \((x_i)\) is weakly null
and \(q \in \mathbb{N}\) was arbitrary.
Concluding, there exists \(L \in [M_n]\) such that 
\([L] \subset {\Delta}_{F}\).
By repeating the previous argument successively over all possible
subsets of \(T_{n+1}\), we obtain
\(M_{n+1} \in [M_n]\) such that \([M_{n+1}] \subset {\Delta}_{F}\),
for all \(F \subset T_{n+1}\).
This completes the inductive construction.
Set \(M=(m_i)\). Let \(L \in [M]\) and \(k \in \mathbb{N}\)
and suppose that there exists \(x^* \in B_{X^*}\) 
which is \(\xi\)-good
for \((\lambda , L , k) \). Now let
\(m_n = \max supp {\xi}_{k}^{L}\)
and \(F = \bigl \{ i < n : m_i \notin 
         \bigcup_{j=1}^{k} supp {\xi}_{j}^{L}
          \bigr \} \).
Our construction yields that \((m_i)_{i\in T_n}\) and 
\(\{ m_i : i \geq n+1 \}\), 
satisfy \((F- E_n)\). We also have, by stability, that 
\[{\xi}_{j}^{L} = {\xi}_{j}^{R_n},
\text{ for all } j \leq k ,  \]
where \(R_n = \{ m_i : i \in T_n \setminus F\} \bigcup
              \{ m_i : i > n+1\}\)
and therefore there exists \(y^* \in B_{X^*}\) 
which is \(\xi\)-good
for \((\lambda , L , k) \) and such that
  \[\sum_{i \in F} |y^*(x_{m_i})|< 
\sum_{i=0}^{n}{\epsilon}_i < \epsilon . \]
But \((x_i)\) is bimonotone and thus we can assume that 
\(y^*(x_i) = 0\) for all \(i > m_n\). Hence,
 \(\sum_{i \in M \setminus \bigcup_{j=1}^{k}supp {\xi}_{j}^{L}}
   |y^*(x_i)| < \epsilon \), as desired.
\end{proof}

\begin{proof}[Proof of Theorem \ref{T:4}]
Let \(\delta > 0\). Our goal is to find \(M \in [\mathbb{N}]\)
and a constant \(C(\delta) > 0\) such that if \(L \in [M]\),
\(n \in \mathbb{N}\), and \((a_i)_{i=1}^{n}\) are scalars in
\([-1,1]\), then
 \[\biggl \| \sum_{i \in F}a_i {\xi}_{i}^{L} \cdot s \biggr \| \leq 
C(\delta) \max \biggl \{ \delta , \, \biggl \| \sum_{i=1}^{n} 
a_i {\xi}_{i}^{L} \cdot s \biggr \| \biggr \},\] 
for all \( F \subset \{ i \leq n : \, |a_i| \geq \delta \} \).
If this is accomplished, then a simple diagonalization argument
yields \(M \in [\mathbb{N}]\) which works for all \(\delta > 0\).
We let \({\lambda}_k = 1 + \frac{k+1}{\delta}\), for all
\(k \in \mathbb{N} \bigcup \{0\}\).
Inductively we construct a decreasing sequence \((M_k)_{k=0}^{\infty}\)
of infinite subsets of \(\mathbb{N}\) such that for all
\(k \in \mathbb{N} \bigcup \{0\}\), \(M_k\) satisfies the
conclusion of Lemma \ref{T:32} for
``$\lambda$'' $= {\lambda}_k$ and ``$\epsilon$'' $= \delta$.
Next choose \(m_1 < m_2 < \cdots \) so that \(m_i \in M_i\),
for all \( i \in \mathbb{N}\). Let \(M=(m_i)\).
We shall show that \(M\) is the desired. To this end, let
\(L \in [M]\) and scalars \((a_i)_{i=1}^{n}\) in \([-1,1]\).
Let also \(F \subset \{1, \cdots ,  n \}\) such that
\(|a_i| \geq \delta \), for all \(i \in F\).
Choose \(k \in \mathbb{N} \bigcup \{0\}\) so that
\[ k \leq \biggl \| \sum_{i=1}^{n}
a_i {\xi}_{i}^{L} \cdot s \biggr \|
< k + 1\]
and set \(F_k = \{ i \in F : i \geq k \} \).
We claim that 
\[ \biggl \| \sum_{i \in F_k}a_i {\xi}_{i}^{L} \cdot s \biggr \|
  \leq 4 {\lambda}_k . \]
Assume this is not the case and choose \(x^* \in B_{X^*}\)
such that \[\bigl | \sum_{i \in F_k}a_i x^*({\xi}_{i}^{L} \cdot s) 
\bigr | > 4 {\lambda}_k . \]
We can further choose \(G_k \subset F_k\) such that 
  \[ \biggl | \sum_{i \in G_k}a_i x^*({\xi}_{i}^{L} \cdot s) \biggr |
 > {\lambda}_k , \]
  \[\text{ the scalars }\, (a_i)_{i \in G_k} \text{ are all of 
    the same sign, and,  } \]
  \[x^*({\xi}_{i}^{L} \cdot s) \geq 0 , \, \text{ for all }
  i \in G_k , (\text{ by replacing } x^* \text{ by }
  -x^* , \text{ if necessary }). \]
Observe that \(supp {\xi}_{i}^{L} \subset M_k \),
when \(k \leq i \leq n\), and thus by stability, there exists
\(y^* \in B_{X^*}\) such that
 \[ y^*({\xi}_{i}^{L} \cdot s) \geq 0 , \, \text{ for all }
  i \in G_k , \]
  \[ \sum_{i \in G_k} y^*({\xi}_{i}^{L} \cdot s) > {\lambda}_k
  \, \text{ and } \,
\sum_{i \in M_k  \setminus \bigcup_{j \in G_k} supp {\xi}_{j}^{L}}
   |y^*(x_i)| < \delta . \]
But now,
\begin{align}
  \biggl | \sum_{i=k}^{n} a_i y^*({\xi}_{i}^{L} \cdot s) \biggr | 
  &\geq  
  \biggl | \sum_{i \in G_k} a_i y^*({\xi}_{i}^{L} \cdot s) \biggr |
  - \sum_{ i \in \{k, \cdots , n\} \setminus G_k}
    |y^*({\xi}_{i}^{L} \cdot s)| \notag \\
  &> \sum_{i \in G_k} |a_i|  y^*({\xi}_{i}^{L} \cdot s)
  - \sum_{i \in M_k \setminus \bigcup_{j \in G_k} supp {\xi}_{j}^{L}}
    |y^*(x_i)| \notag \\
  &> {\delta}{{\lambda}_k} - \delta = k + 1. \notag
\end{align}
Thus, \(\|\sum_{i=k}^{n} a_i {\xi}_{i}^{L} \cdot s \| > k + 1\)
and since \((x_n)\) is bimonotone, we also have that
\(\|\sum_{i=1}^{n} a_i {\xi}_{i}^{L} \cdot s \| > k + 1\)
which is a contradiction. Therefore, our claim holds and hence
\[\|\sum_{i \in F} a_i {\xi}_{i}^{L} \cdot s \| \leq
\begin{cases}
4(1 + \frac{1}{\delta}), &\text{ if } 
\|\sum_{i=1}^{n} a_i {\xi}_{i}^{L} \cdot s \| < 1; \\
k(5 + \frac{8}{\delta}), &\text{ if } k \leq 
\|\sum_{i=1}^{n} a_i {\xi}_{i}^{L} \cdot s \| < k + 1, \, k \in \mathbb{N}.
\end{cases}
\]    
Concluding,
\[\biggl \|\sum_{i \in F} a_i {\xi}_{i}^{L} \cdot s \biggr \| \leq
  C(\delta)
\max \biggl \{ \delta , \, \biggl \| \sum_{i=1}^{n} 
a_i {\xi}_{i}^{L} \cdot s \biggr \| \biggr \},\] 
where, \(C(\delta)=\frac{1}{\delta} (5 + \frac{8}{\delta})\).
The proof of Theorem \ref{T:4} is now complete.
\end{proof}
\begin{remark}
A refinement of the proof of Lemma \ref{T:32}, yields
that given \(\theta > 0\), there exists \(M \in [\mathbb{N}]\)
satisfying the following property:
If \(L \in [M]\), \(k \in \mathbb{N}\) 
and there exists \(x^* \in B_{X^*}\) which is \(\xi\)-good for
\((\lambda , L , k)\), then there exists \(y^* \in B_{X^*}\)
with \(y^*({\xi}_{i}^{L} \cdot s) \geq 0\) for \(i \leq k\)
and such that
\begin{align}
&\sum_{i=1}^{k} y^*( {\xi}_{i}^{L} \cdot s) \geq
 (1 - \theta )\sum_{i=1}^{k} x^*( {\xi}_{i}^{L} \cdot s)
\text{ and, } \notag \\
&\sum_{n \in M \setminus \bigcup_{i=1}^{k}supp {\xi}_{i}^{L}}
   |y^*(x_n)| < \epsilon \sum_{i=1}^{k} y^*( {\xi}_{i}^{L} \cdot s).
\notag 
\end{align}
This property in turn implies that the constant \(C(\delta)\)
found in the proof of Theorem \ref{T:4}, is actually
proportional to \(\frac{1}{\delta}\).
\end{remark}
Elton's aforementioned theorem follows immediately
if we let \(\xi = 0\) in the statement of Theorem \ref{T:4}.
As a corollary to Elton's theorem one obtains the following
dichotomy \cite{e}, \cite{o1}:
\begin{Cor} \label{T:33}
A normalized weakly null sequence \((x_n)\) in a Banach space, either
has a subsequence equivalent to the unit vector basis of
\(c_0\), or, a Schauder basic subsequence \((x_{m_n})\) 
which is semi-boundedly complete.
If the latter alternative occurs, then the sequence of the biorthogonal
functionals to \((x_{m_n})\), converges weakly to zero in the dual of
the closed linear span of \((x_{m_n})\).
\end{Cor}

Our next result is the dual form of Lemma \ref{T:32}.
\begin{Lem} \label{T:34}
Let \(s = (x_n)\) be a weakly null sequence in the Banach space \(X\)
and \(\xi < {\omega}_1\).
Given \(\epsilon > 0\), \(\delta > 0\), \(N \in [\mathbb{N}]\) there exists 
\(M \in [N]\) satisfying the following
property: If \(k \in \mathbb{N}\), \(L \in [M]\),  and there exists
\(x^* \in B_{X^*}\) so that
\(x^*({\xi}_{n}^{L} \cdot s) \geq \delta \),
for all \( n \leq k\), then there exists
\(y^* \in B_{X^*}\) so that
\[y^*({\xi}_{n}^{L} \cdot s) \geq \delta , 
  \text{ for all } n \leq k , \, \text{ and } , \,
\sum_{n \in M \setminus \bigcup_{i=1}^{k}supp {\xi}_{i}^{L}}
   |y^*(x_n)| < \epsilon. \]
\end{Lem}
\begin{proof}
The proof is similar to that of Lemma \ref{T:32}.
Let us now say that the functional \(x^* \in B_{X^*}\)
is good for \((\delta , L , k)\), if
\(x^*({\xi}_{n}^{L} \cdot s) \geq \delta \), for all
\(n \leq k\).
Next choose \(({\epsilon}_i)_{i=0}^{\infty}\), 
a sequence of positive scalars
such that \(\sum_{i=0}^{\infty} {\epsilon}_i < \epsilon\).
Using the same notation and terminology as in Lemma \ref{T:32},
let \(n \in \mathbb{N} \bigcup \{0\}\) and 
\(F \subset T_n\). 
The \(n\)-tuple \((m_i)_{i \in T_n}\)
and the infinite subset \(L\) of \(\mathbb{N}\),
(\(L=(l_i)\)), are said to
satisfy property \((F-E_n^*)\), provided that
\(m_n < l_1\), if \(n \geq 1\), and the following statement holds:

If \(k \in \mathbb{N}\),
and there exists \(x^* \in B_{X^*}\) which is good for
\((\lambda , \{m_j : j \in T_n \setminus F\} \bigcup
\{l_j : j \geq 2 \} , k)\), 
then there exists
\(y^* \in B_{X^*}\) which is good for
\((\lambda , \{m_j : j \in T_n \setminus F\}  \bigcup
\{l_j : j \geq 2 \} , k)\)
and such that
\[ \sum_{j \in F} |y^*(x_{m_j})| + |y^*(x_{l_1})|
< \sum_{i=0}^{n}{\epsilon}_i .\] 
Let us also say that \((m_i)_{i \in T_n}\) and \(L\)
satisfy property \((E_n^*)\), if they satisfy property
\((F-E_n^*)\), for every \(F \subset T_n \). 

We shall inductively construct an increasing sequence
\((m_n)_{n=1}^{\infty}\) 
of elements of \(N\), and a decreasing sequence
\((M_n)_{n=0}^{\infty}\) of infinite subsets of \(N\) 
with \(m_n \in M_{n-1}\), if \(n \geq 1\), and
so that 
for every \(n \in \mathbb{N}\bigcup \{0\}\), 
if \(L \in [M_n]\), then
\((m_i)_{i \in T_n}\) and \(L\) satisfy property \((E_n^*)\).

The first inductive step is similar to the general one and
therefore we shall not explicit it. Assume that 
\((m_i)_{i \in T_n}\) and 
\(M_0 \supset \cdots \supset M_n \), infinite subsets of \(N\),
with \(m_i \in M_{i-1}\) for \(1 \leq i \leq n\)
have been constructed so that if \(i \leq n\) and 
\(L \in [M_i]\), then  \((m_j)_{j \in T_i}\) and \(L\) 
satisfy property \((E_i^*)\). 
Let \(m_{n+1}= \min M_n\). Fix \(F \subset T_{n+1}\)
and define
\[ {\Delta}_F = \{ L \in [M_n], L=(l_i),: (m_i)_{i \in T_{n+1}}
   \text{ and } L \text{ satisfy } (F-E_{n+1}^*) \}.\]
Clearly, \({\Delta}_F\) is closed in \([M_n]\) and therefore
Ramsey. Arguing as in the proof of Lemma \ref{T:32}, we
obtain \(M_{n+1} \in [M_n]\) such that
\([M_{n+1}] \subset {\Delta}_F\), for every \( F \subset T_{n+1}\).
Indeed, we need only modify the definition of \(k_{j_0}\) in the argument
of Lemma \ref{T:32}. We alternatively set
\(k_{j_0}= \max \{ k_j : j \leq q\}\)
 and observe that if \(y^*\) is good for \((\delta , R , k_{j_0} )\),
then \(y^*\) is also good for \((\delta , R , k_j )\), for
every \( j \leq q\).
The argument of Lemma \ref{T:32} is now carried over unaltered
yielding the proof of Lemma \ref{T:34}.
\end{proof}

In the proof of the Schreier unconditionality theorem,
\cite{mr}, \cite{o2},
we shall make use
of the following:
\begin{Lem} \label{T:35}
Let \((x_n)\) be a weakly null sequence in the Banach space \(X\).
Let also \(\epsilon > 0\), \(\delta > 0\) and \(k \in \mathbb{N}\).
There exists \(Q \in [\mathbb{N}]\), \(Q=(q_n)\), such that for
every \(x^* \in B_{X^*}\) and \(F \in [\mathbb{N}]^k\)
there exists \(y^* \in B_{X^*}\) satisfying
\[ \sum_{i \in F} |x^*(x_{q_i}) - y^*(x_{q_i})| < \delta , \,
   \text{ and } , \,
   \sum_{i \notin F} |y^*(x_{q_i})| < \epsilon . \]
\end{Lem}
\begin{proof}
Again, the proof is much similar to that of Lemma \ref{T:32}.
Choose first \(({\epsilon}_i)_{i=0}^{\infty}\), 
a sequence of positive scalars
such that \(\sum_{i=0}^{\infty} {\epsilon}_i < \epsilon\).
Let \(\Vec{\lambda}= ({\lambda}_1, \cdots , {\lambda}_k)\)
be an element in \([-1,1]^k\). We shall say that the functional
\(x^* \in B_{X^*}\) is \(\Vec{\lambda}\)-good for 
\(L\), where \(L=(l_i)\) is an infinite subset of
\(\mathbb{N}\), if
\(\sum_{i=1}^{k}|x^*(x_{l_i}) -{\lambda}_i | < \delta\). 
Using the same notation and terminology as in Lemma \ref{T:32},
let \(n \in \mathbb{N} \bigcup \{0\}\) and 
\(F \subset T_n\). 
The \(n\)-tuple \((m_i)_{i \in T_n}\)
and the infinite subset \(L\) of \(\mathbb{N}\),
(\(L=(l_i)\)), are said to
satisfy property \((F-O_n)\), provided that
\(m_n < l_1\), if \(n \geq 1\), and the following statement holds:
 
If  there exists \(x^* \in B_{X^*}\) which is \(\Vec{\lambda}\)- good for
\(\{m_j : j \in T_n \setminus F\} \bigcup
\{l_j : j \geq 2 \}\), 
then there exists
\(y^* \in B_{X^*}\) which is \(\Vec{\lambda}\)-good for
\(\{m_j : j \in T_n \setminus F\}  \bigcup
\{l_j : j \geq 2 \}\)
and such that
\[ \sum_{j \in F} |y^*(x_{m_j})| + |y^*(x_{l_1})|
< \sum_{i=0}^{n}{\epsilon}_i .\] 
Let us also say that \((m_i)_{i \in T_n}\) and \(L\)
satisfy property \((O_n)\), if they satisfy property
\((F-O_n)\), for every \(F \subset T_n \). 
We shall inductively construct an increasing sequence
\((m_n)_{n=1}^{\infty}\) 
of elements of \(\mathbb{N}\), and a decreasing sequence
\((M_n)_{n=0}^{\infty}\) of infinite subsets of \(\mathbb{N}\) 
with \(m_n \in M_{n-1}\), if \(n \geq 1\), 
so that 
for every \(n \in \mathbb{N}\bigcup \{0\}\), 
if \(L \in [M_n]\), then
\((m_i)_{i \in T_n}\) and \(L\) satisfy property \((O_n)\).

The first inductive step is similar to the general one and
therefore we shall not explicit it. Assume that 
\((m_i)_{i \in T_n}\) and 
\(M_0 \supset \cdots \supset M_n \), infinite subsets of \(N\),
with \(m_i \in M_{i-1}\) for \(1 \leq i \leq n\)
have been constructed so that if \(i \leq n\) and 
\(L \in [M_i]\), then  \((m_j)_{j \in T_i}\) and \(L\) 
satisfy property \((O_i)\). 
Let \(m_{n+1}= \min M_n\). Fix \(F \subset T_{n+1}\)
and define
\[ {\Delta}_F = \{ L \in [M_n], L=(l_i),: (m_i)_{i \in T_{n+1}}
   \text{ and } L \text{ satisfy } (F-O_{n+1}) \}.\]
Clearly, \({\Delta}_F\) is closed in \([M_n]\) and therefore
Ramsey. Arguing as in the proofs of Lemmas \ref{T:32} and \ref{T:34}       
we obtain \(M_{n+1} \in [M_n]\) such that
\([M_{n+1}] \subset {\Delta}_F\), for every \( F \subset T_{n+1}\).
The inductive construction is now complete and we set
\(M=(m_i)\). It follows, by our construction, that if 
\(F \in [\mathbb{N}]^k\) and
\(x^* \in B_{X^*}\) is \(\Vec{\lambda}\)-good for
\(\{m_j : j \in F \} \bigcup \{ m_j : j > \max F\}\),
then there exists
\(y^* \in B_{X^*}\) \(\Vec{\lambda}\)-good for
\(\{m_j : j \in F \} \bigcup \{ m_j : j > \max F\}\)
and such that 
\(\sum_{j \notin F}|y^*(x_{m_j})| < \epsilon\).
Let us then say that \(M\) works for \(\Vec{\lambda}\).
Finally, let \(\mathcal{E}\) be a finite \(\delta\)-net
in \([-1,1]^k\), and choose \(Q \in [\mathbb{N}]\)
which works for every \(\Vec{\lambda} \in \mathcal{E}\).
It is easily verified that \(Q\) is the desired.
\end{proof}

Lemma \ref{T:34} provides an alternative proof of the
fact that every normalized weakly null sequence admits a 
convexly unconditional subsequence.
\begin{proof}[Proof of Theorem \ref{T:5}.]
Let \(\delta > 0 \). It is enough to find \(M \in [\mathbb{N}]\)
, \(M=(m_i)\), so that if \(F \in [\mathbb{N}]^{< \infty}\)
and \(({\lambda}_i)_{i \in F}\) are scalars satisfying
\(\|\sum_{i \in F}{\lambda}_i x_{m_i} \| > \delta \)
and \(\sum_{i \in F}|{\lambda}_i| \leq 1\),
then, \(\|\sum_{i=1}^{\infty} a_i x_{m_i} \| > \frac{{\delta}^2}{32}\),
for all choices of scalars \((a_i)_{i=1}^{\infty} \subset c_{00}\)
with \(\max_{i}|a_i| \leq 1\) and 
such that \(|a_i|=|{\lambda}_i|\), for all \(i \in F\).
Once this is accomplished, then a simple diagonalization argument
yields \(M \in [\mathbb{N}]\) which works for all \(\delta > 0\).
To this end, let \(M=(m_i)\) be the infinite subset of 
\(\mathbb{N}\) resulting from Lemma \ref{T:34} applied
on the sequence \((x_n)\) for \(\xi=0\), 
``\(\delta\)''\( = \frac{\delta}{2}\) 
and ``\(\epsilon\)''\(=\frac{{\delta}^2}{32}\).
Let now \(F \in [\mathbb{N}]^{< \infty}\) and 
\(({\lambda}_i)_{i \in F}\) be scalars such that
\(\|\sum_{i \in F}{\lambda}_i x_{m_i} \| > \delta \) and
\(\sum_{i \in F}|{\lambda}_i| \leq 1\).
Choose \(x^* \in B_{X^*}\) such that
\(|\sum_{i \in F}{\lambda}_i x^*(x_{m_i}) | > \delta \) and set
\[G = \{ i \in F : | x^*(x_{m_i}) | \geq \frac{\delta}{2} \}. \]
Clearly, \(|\sum_{i \in G}{\lambda}_i x^*(x_{m_i}) | >\frac{\delta}{2} \).
Next, let \((a_i)_{i=1}^{\infty} \subset c_{00}\)
with \(\max_{i}|a_i| \leq 1\) and 
\(|a_i|=|{\lambda}_i|\), for all \(i \in F\).
By splitting \(G\) into four sets in the obvious manner, we find
\(H \subset G\) such that
  \[ \biggl | \sum_{i \in H}{\lambda}_i 
   x^*(x_{m_i}) \biggr |
 > \frac{\delta}{8}, \]
  \[\text{ the scalars }\, ({\lambda}_i )_{i \in H} \text{ are all of 
    the same sign, and,  } \]
  \[x^*(x_{m_i}) \geq \frac{\delta}{2} \, \text{ for all }
  i \in H , (\text{ by replacing } x^* \text{ by }
  -x^* , \text{ if necessary }). \]
Now choose \(y^* \in B_{X^*}\) with
\(y^*(x_{m_i}) \geq \frac{\delta}{2}\),
for all \(i \in H\), and such that
\(\sum_{i \notin H}|y^*(x_{m_i})| < \frac{{\delta}^2}{32}\).
It follows that 
\begin{align}
\biggl \|\sum_{i=1}^{\infty} a_i x_{m_i} \biggr \|  &\geq
 \biggl | \sum_{i \in H}a_i 
   y^*(x_{m_i}) \biggr | - \sum_{i \notin H}|y^*(x_{m_i})| \notag \\
  &> \sum_{i \in H}|{\lambda}_i| \frac{\delta}{2} - \frac{{\delta}^2}{32}
   > \frac{{\delta}^2}{32}. \notag
\end{align}
\end{proof}

\begin{Cor} \label{T:36} \hfill
For a normalized weakly null sequence \(s=(x_n)\) and 
\(\xi < {\omega}_1\), the following are equivalent:
\begin{enumerate}
\item There exists \(M \in [\mathbb{N}]\), \(M=(m_i)\), 
so that \((x_{m_i})\) is an \(l_{1}^{\xi}\) spreading model.
\item There exist \(N \in [\mathbb{N}]\) and \(\delta > 0\)
such that \(S_{\xi}(N) \subset {\mathfrak{F}}_{\delta}\).
\end{enumerate}
\end{Cor}
\begin{proof}
Suppose 1. holds and let \(C > 0\) such that
\[\|\sum_{i \in F}{\lambda}_i x_{m_i} \| \geq 
 C \sum_{i \in F}|{\lambda}_i|
\text{ for all } F  \in S_{\xi},\, \text{ and scalars }
({\lambda}_i)_{i \in F}.\]
Let \(t\) denote the sequence \((y_i)\), where
\(y_i = x_{m_i}\), for all \(i \in \mathbb{N}\).
Then, \(\|{\xi}_{1}^{L} \cdot t\| \geq C\),
for all \(L \in [\mathbb{N}]\). It follows that for
\(\delta = \frac{C}{2}\), the hereditary family
\[ \{F \in [\mathbb{N}]^{< \infty} : \,
     M(F) \in {\mathfrak{F}}_{\delta}\},\]
is \((1, \xi, \mathbb{N}, \delta)\) large. Corollary
\ref{T:217} now yields \(N \in [M]\) so that
\(S_{\xi}(N) \subset {\mathfrak{F}}_{\delta}\)
and thus 2. holds.

Assume now that 2. holds and choose \(M = (m_i) \in [N]\)
as in the proof of Theorem \ref{T:5}, applied on the
sequence \((x_n)_{n \in N}\) for 
``\(\delta\)''\(= \frac{\delta}{4}\).
Let \(F \in S_{\xi}\) and scalars \(({\lambda}_i)_{i \in F}\)
such that \(\sum_{i \in F}|{\lambda}_i|=1\).
We claim that \(\|\sum_{i \in F}{\lambda}_i x_{m_i} \|
\geq \frac{{\delta}^2}{512}\), which evidently yields 1.
Indeed, by our assumption, there exists
\(x^* \in B_{X^*}\) such that \(|x^*(x_{m_i})| \geq \delta\),
for all \(i \in F\). Next choose \( G \subset F\) such that
\[ \sum_{i \in G}|{\lambda}_i| \geq \frac{1}{4}, \]
\[\text{ the scalars } ({\lambda}_i)_{i \in G} \text{ are all
of the same sign, and, } \]
\[x^*(x_{m_i}) \geq \delta , \, \text{ for all } i \in G, \,
(\text{ by replacing } x^* \text{ by } -x^* \text{ if necessary }). \]
Therefore, \(\|\sum_{i \in G}{\lambda}_i x_{m_i} \| \geq
  |\sum_{i \in G}{\lambda}_i x^*(x_{m_i}) | 
\geq \frac{\delta}{4}\), and hence applying Theorem \ref{T:5}
we obtain that \(\|\sum_{i \in F}{\lambda}_i x_{m_i}\| \geq 
\frac{1}{32}(\frac{\delta}{4})^2=\frac{{\delta}^2}{512}\),
as claimed.
\end{proof}

\begin{proof}[Proof of Theorem \ref{T:3}.]
Assume that 2. does not hold. Let \(M\in [\mathbb{N}]\)
and \(\epsilon > 0\). It is easily seen that the set
\[{\mathcal{A}}_{\epsilon} = \{ L \in [M]: 
\|{\xi}_{1}^{L} \cdot s \| < \epsilon \} \]
where \(s=(x_n)\), is closed in \([M]\) and therefore Ramsey.
If it were the case that \([P] \bigcap {\mathcal{A}}_{\epsilon}=
\emptyset\), for some \(P \in [M]\), then the family
\({\mathfrak{F}}_{\frac{\epsilon}{2}}\)
would be \((1, \xi, P ,\frac{\epsilon}{2})\) large, and hence,
by Corollaries \ref{T:217} and \ref{T:36},
\((x_n)_{n \in L}\) would in turn be an \(\ell_1^{\xi}\)
spreading model, for some \(L \in [P]\) contradicting
our assumption. It follows now that we can construct
\((M_n)\), a decreasing sequence of infinite subsets of \(\mathbb{N}\),
such that for all \( n \in \mathbb{N}\),
\[\|{\xi}_{1}^{L} \cdot s \| < \frac{1}{n},
\text{ for all } L \in [M_n]. \]
Let now \(N\) be any infinite subset of \(M\) almost contained
in each \(M_n\), and it is easy to verify that \((x_n)\)
is \((L, \xi)\) convergent, for all \(L \in [N]\).

In order to show that 1. and 2. are mutually exclusive, assume that 
\((x_{m_i})\) is an \(\ell_1^{\xi}\)
spreading model with constant \(C \).
We can choose \(N \in [M]\) such that 
\begin{align}
F \setminus \{\min F\} \in S_{\xi}(M), \, \text{ for all }
   F \in  S_{\xi}[N], \, &(\text{ by Lemma } \ref{T:24}), \notag \\ 
\|{\xi}_{1}^{L}\|_0 < \epsilon , \, \text{ for all }
    L \in [N], \, &(\text{ by Proposition } \ref{T:215}). \notag 
\end{align}
where \(\epsilon > 0\) is chosen so that 
\(\epsilon < \frac{C}{1+C}\).
It follows now that for every \(L \in [N]\),
\[ \|{\xi}_{n}^{L} \cdot s \| > C - (1+C){\epsilon} > 0,
\text{ for all } n \in \mathbb{N}, \]
and thus \((x_n)\) is not \((L,\xi)\) convergent.
Hence 1. does not hold.
\end{proof}
An immediate consequence of Theorem \ref{T:3} is Corollary \ref{C:1}.
\begin{proof}[Proof of Corollary \ref{C:1}.]
Suppose that \(s=(x_n)\) is an \(\ell_1^{\xi}\)
spreading model with constant \(C\). It follows that
\(\|{\xi}_n^L \cdot s\| \geq C \), for all \(n \in \mathbb{N}\)
and \(L \in [\mathbb{N}]\). Next, choose according
to Theorem \ref{T:3}, \(N \in [\mathbb{N}]\)
so that \((x_n)\) is \((L, \xi +1 )\) convergent, for
every \(L \in [N]\). Evidently, 
\(({\xi}_n^L \cdot s)\) is Cesaro summable for all \(L \in [N]\).
\end{proof}

We continue our applications of Theorem \ref{T:5} with the
\begin{proof}[Proof of Theorem \ref{T:6}.]
Let \(C\) be the constant of the \(\ell_1^{\xi}\)
spreading model \((x_n)\).
Let now \(M \in [\mathbb{N}]\), \(M=(m_i)\) satisfying the
conclusion of Theorem \ref{T:5} for the sequence \((x_n)\)
and ``\(\delta\)''\( = \frac{C}{2}\).
We claim that \((x_{m_i})\) is \(S_{\xi}\) unconditional.
Indeed, let \( n \in \mathbb{N}\) and scalars 
\((a_i)_{i=1}^{n}\) be given. Let also
\(F \subset \{1, \cdots , n\}\), \(F \in S_{\xi}\),
such that \(\|\sum_{i \in F}a_i x_{m_i} \|=1\).
It follows that \(C \sum_{i \in F} |a_i| \leq 1\)
and \(\|\sum_{i \in F} C a_i x_{m_i} \| >\frac{C}{2} \).
If there exists \(j \leq n\) so that
\(C |a_j| > 1\), then, since \((x_n)\) is bimonotone,
we obtain that \(\|\sum_{i=1}^n a_i x_{m_i} \| > \frac{1}{C}\).
So assuming that \(C |a_i| \leq 1\), for all \(i \leq n\),
we obtain through Theorem \ref{T:5} that
\(\|\sum_{i=1}^n a_i x_{m_i} \| > \frac{C}{128}\).
Hence, 
\[ \biggl \|\sum_{i \in F}a_i x_{m_i} \biggr \| \leq \frac{128}{C}
   \biggl \|\sum_{i=1}^n a_i x_{m_i} \biggr \|,\]
for every \(F \in S_{\xi}\) and all choices of scalars \((a_i)_{i=1}^{n}\).
\end{proof}

We also obtain the result on Schreier unconditionality,
\cite{mr}, \cite{o2}.
\begin{Cor}
Let \((x_n)\) be a normalized weakly null sequence in \(X\) and
\(\epsilon > 0\). There exists a subsequence \((x_{m_i})\)
of \((x_n)\) which is \(2 + \epsilon \) \(S_1\) unconditional.
\end{Cor}
\begin{proof}
Choose first \(\theta > 0\) such that 
\((1 + \theta)(2 + \theta) < 2 + \epsilon\).
By passing to a subsequence, if necessary, 
we can assume that \((x_n)\) is Schauder basic with basis constant
\(1 + \theta\). We first show that for every \(k \in \mathbb{N}\)
and \(N \in [\mathbb{N}]\), there exists \(L \in [N]\),
\(L=(l_i)\) so that
\[ \biggl \| \sum_{i \in F}a_i x_{l_i} \biggr \| \leq
 (1 + \theta) \biggl \| \sum_{i=1}^{\infty} a_i x_{l_i} \biggr \|, \]
for every \(F \in [\mathbb{N}]^k\) and all choices of scalars
\((a_i)\) in \(c_{00}\).
Indeed, apply Lemma \ref{T:35} to the sequence \((x_n)_{n \in N}\)
to obtain \(L \in [N]\), \(L=(l_i)\), satisfying the conclusion
of that lemma for \( k\), 
``\(\delta\)''\( = \frac{\theta}{6}\), and
``\(\epsilon\)''\( = \frac{\theta}{6}\).
Let now \((a_i) \subset c_{00}\) such that
\(\| \sum_{i=1}^{\infty} a_i x_{l_i}\|=1\).
Let also \(F \in [\mathbb{N}]^k\) and choose
\(x^* \in B_{X^*}\) such that
\[\biggl \| \sum_{i \in F}a_i x_{l_i} \biggr \|=
  \biggl | \sum_{i \in F}a_i x^*(x_{l_i}) \biggr |.\]
Then choose \(y^* \in B_{X^*}\) such that
\[\sum_{i \in F}|x^*(x_{l_i}) - y^*(x_{l_i}) | < 
  \frac{\theta}{6}, \, \text{ and }, \, 
  \sum_{i \notin F}|y^*(x_{l_i}) | <\frac{\theta}{6}. \]
We now have the following estimate
\begin{align}
 1= &\biggl \| \sum_{i=1}^{\infty} a_i x_{l_i} \biggr \| \geq
 \biggl | \sum_{i=1}^{\infty} a_i y^*(x_{l_i}) \biggr | \notag \\
 &\geq \biggl | \sum_{i \in F}a_i x^*(x_{l_i}) \biggr | - 
  \sum_{i \in F} |a_i| |x^*(x_{l_i}) - y^*(x_{l_i}) | -
  \sum_{i \notin F} |a_i| |y^*(x_{l_i}) | \notag \\
 &\geq \biggl \| \sum_{i \in F}a_i x_{l_i} \biggr \| -
   3 \frac{\theta}{6} - 3 \frac{\theta}{6}, \notag
\end{align}
and thus,
\( \| \sum_{i \in F}a_i x_{l_i}\| \leq (1 + \theta)
\| \sum_{i=1}^{\infty} a_i x_{l_i}\|\), as desired.
We can now choose a decreasing sequence \((M_k)\)
of infinite subsets of \(\mathbb{N}\) such that
for every \(k \in \mathbb{N}\),
\[  \biggl \| \sum_{m \in F}a_m x_m \biggr \| \leq
 (1 + \theta) \biggl \| \sum_{m \in M_k}a_m x_m \biggr \|, \] 
for all \( F \in [M_k]^j\) and  \( j \leq k \),
and all choices of scalars \((a_m)_{m \in M_k} \subset c_{00}\).
Finally choose \(m_1 < m_2 < \cdots \) with
\(m_i \in M_i\), for all \(i \in \mathbb{N}\), and set
\(M=(m_i)\). It is easily verified that the subsequence
\((x_{m_i})\) is \(2 + \epsilon\) \(S_1\) unconditional.
\end{proof} 

The final results of this section concern the duality between
\({c_0}^{\xi}\) and \(\ell_1^{\xi}\) spreading models. We first
recall the following
\begin{Def}
A sequence \((x_n)\) in a Banach space is called a \({c_0}^{\xi}\)
spreading model, \(1 \leq \xi < {\omega}_1\), if there exists a
constant \(C > 0\) such that
\[ \biggl \|\sum_{i \in F} a_i x_i \biggr\| \leq
C \max_{i \in F}|a_i|,\]
for every \(F \in S_{\xi}\) and all choices of scalars
\((a_i)_{i \in F}\).
\end{Def}

\begin{notation}
If \(M \in [\mathbb{N}]\) and \((x_n)\) is a sequence in \(X\)
then we let \(X_M\) denote the closed linear span of the subsequence
\((x_n)_{n \in M}\).
\end{notation}
The duality between \({c_0}^{\xi}\) and \(\ell_1^{\xi}\) spreading models
is described in the following
\begin{Thm} \label{T:39}
Let \((x_n)\) and \((x_n^*)\) be normalized weakly null sequences
in \(X\) and \(X^*\) respectively. Assume that for some \(\epsilon > 0\)
we have that \(\inf_n |x_n^*(x_n)| \geq \epsilon\). Let 
\(1 \leq \xi < {\omega}_1\). The following are equivalent:
\begin{enumerate}
\item For every \(M \in [\mathbb{N}]\), there exists \(L \in [M]\)
such that \((x_n^*|X_M)_{n \in L}\) is an \(\ell_1^{\xi}\)
spreading model.
\item For every \(M \in [\mathbb{N}]\), there exists \(L \in [M]\)
such that \((x_n)_{n \in L}\) is a \({c_0}^{\xi}\) spreading model.
\end{enumerate}
\end{Thm}
\begin{proof}
Once again, we assume our sequence \((x_n)\) is bimonotone.
We can also assume, without loss of generality, that
\[ \sum_{i=1}^{\infty} \sum_{j \ne i} |x_i^*(x_j)| < \frac{\epsilon}{2}.\]
Furthermore, we shall assume that \((x_n)\) satisfies the 
conclusion of Theorem \ref{T:4} for \(\xi=0\) and \(M = \mathbb{N}\).
That is, for all \(\delta > 0\). there exists a constant
\(C(\delta) > 0\) such that for every \(n \in \mathbb{N}\)
and all scalars \((a_i)_{i=1}^n\) in \([-1,1]\), if
\(F \subset \{1, \cdots , n\}\) and \(|a_i| \geq \delta\) for all
\(i \in F\), then
\[\biggl \| \sum_{i \in F} a_i x_i \biggr \| \leq C(\delta)
 \biggl \| \sum_{i=1}^n a_i x_i \biggr \|. \]

Assume first that 2. holds and choose \(L \in [M]\),
\(L=(l_i)\), and \(C > 0\) so that
\[ \biggl \|\sum_{i \in F} a_i x_{l_i} \biggr\| \leq
C \max_{i \in F}|a_i|,\]
for every \(F \in S_{\xi}\) and all choices of scalars
\((a_i)_{i \in F}\).
We claim that \((x_{l_i}^*|X_M)\) is an \(\ell_1^{\xi}\)
spreading model. Indeed, let \(F \in S_{\xi}\) and scalars
\((a_i)_{i=1}^{n}\). For each \(i \in F\), let \({\epsilon}_i\)
be the sign of \(a_i x_{l_i}^*(x_{l_i})\).
Of course, \(\|\sum_{i \in F}{\epsilon}_i x_{l_i}\| 
\leq C \). Therefore,
\begin{align}
C \biggl \| \sum_{i \in F} a_i x_{l_i}^*|X_M \biggr \| &\geq
 \biggl | \sum_{i \in F} a_i x_{l_i}^*(\sum_{j \in F}
 {\epsilon}_j x_{l_j} ) \biggr | \notag \\
&\geq \sum_{i \in F}|a_i||x_{l_i}^*(x_{l_i})| -
      \sum_{i \in F}|a_i| \sum_{j \ne i}|x_{l_i}^*(x_{l_j})| \notag \\
&\geq \frac{\epsilon}{2}\sum_{i \in F}|a_i| \notag 
\end{align} 
and thus 1. holds.

Suppose now that 1. holds. Choose, according to Corollary \ref{T:36},
\(L \in [M]\), \(L= (l_i)\) and \(\delta > 0\) such that
\[S_{\xi}(L) \subset \{ F \in [\mathbb{N}]^{< \infty}:\, \exists \,
 x^{**} \in B_{X_{M}^{**}} \, \text{ with }\, 
|x^{**}(x_n^*)| \geq \delta , \, \forall n \in F\}. \]
It follows that 
\[S_{\xi}(L) \subset \{ F \in [\mathbb{N}]^{< \infty}:\, \exists \,
x \in B_{X_M}\, \text{ with }\, |x_n^*(x)| > \frac{\delta}{2},\, 
\forall n \in F\}. \]
We now claim that for every \(x \in X_M\), the sequence
\(\bigl (\frac{x_{m_i}^*(x)}{x_{m_i}^*(x_{m_i})} x_{m_i}\bigr )\),
where \(M=(m_i)\), is strongly bounded.
Indeed, let \(x=\sum_{i=1}^{\infty}c_ix_{m_i} \in B_{X_M}\) and note
that by monotonicity we have that
\(|c_i| \leq 1\), for all \(i \in \mathbb{N}\). 
Then, for all \(k \in \mathbb{N}\),
\begin{align}
\biggl \| \sum_{i=1}^k 
\frac{x_{m_i}^*(x)}{x_{m_i}^*(x_{m_i})} x_{m_i} \biggr \| &\leq
\biggl \| \sum_{i=1}^k c_ix_{m_i} \biggr \| +
\biggl \| \sum_{i=1}^k
\frac{1}{x_{m_i}^*(x_{m_i})} \biggl (\sum_{j \ne i}
 c_jx_{m_i}^*(x_{m_j}) \biggr ) x_{m_i}
\biggr \| \notag \\
&\leq \|x\| + \sum_{i=1}^k \frac{1}{\epsilon}
\sum_{j \ne i}|x_{m_i}^*(x_{m_j})| < \frac{3}{2} \notag
\end{align}
and our claim holds.
Let now \(F \in S_{\xi}\) and choose \(x \in B_{X_M}\)
such that \(|x_{l_i}^{*}(x)| > \frac{\delta}{2}\), for all
\(i \in F\). Since 
\(\bigl |\frac{x_{l_i}^*(x)}{x_{l_i}^*(x_{l_i})} \bigr |
 > \frac{\delta}{2}\), for all \(i \in F\),
our initial assumptions on the sequence \((x_n)\) yield that
\[
\biggl \| \sum_{i \in F} 
\frac{x_{l_i}^*(x)}{x_{l_i}^*(x_{l_i})} x_{l_i} \biggr \|
\leq C \bigl ( \frac{{\delta}{\epsilon}}{2} \bigr)
\biggl \| \sum_{i=1}^n 
\frac{x_{m_i}^*(x)}{x_{m_i}^*(x_{m_i})} x_{m_i} \biggr \|, \]
where \(m_n= l_{\max F}\).
Thus, letting \(b_i^F = \frac{x_{l_i}^*(x)}{x_{l_i}^*(x_{l_i})}\),
for \(i \in F\), we have that
\[
\biggl \| \sum_{i \in F}b_i^F x_{l_i}\biggr \| \leq
\frac{3}{2} C \bigl ( \frac{{\delta}{\epsilon}}{2} \bigr)\]
and \(\frac{\delta}{2} \leq |b_i^F| \leq \frac{1}{\epsilon}\),
for all \(i \in F\) and \(F \in S_{\xi}\).
A standard argument now shows that
\[ \biggl \| \sum_{ i \in F} a_i x_{l_i} \biggr \| \leq
\frac{6}{\delta} C \bigl ( \frac{{\delta}{\epsilon}}{2} \bigr)
\max_{i \in F} |a_i|, \]
for all \(F \in S_{\xi}\) and scalars \((a_i)_{i \in F}\).
Hence, 1. implies 2.
\end{proof}

\begin{Def}
The Banach space \(X\) satisfies the \(\xi\) Dunford-Pettis
property (\(\xi\)-DP), \(1 \leq \xi < {\omega}_1\), if
for every pair of weakly null sequences
\((x_n)\) and \((x_n^*)\) in \(X\) and \(X^*\) respectively,
with \((x_n^*)\) \(\xi\)-convergent, we have that 
\(\lim_n x_n^*(x_n) = 0\).

\(X\) is said to be hereditarily \(\xi\)-DP, if every subspace
of \(X\) satisfies the \(\xi\)-DP.
\end{Def}

\begin{Cor} \label{T:40}
For a Banach space \(X\) and \(1 \leq \xi < {\omega}_1\),
the following are equivalent:
\begin{enumerate}
\item Every normalized weakly null sequence in \(X\) admits 
a subsequence which is a \(c_0^{\xi}\) spreading model.
\item \(X\) is hereditarily \(\xi\)-DP.
\end{enumerate}
\end{Cor}
\begin{proof}
Assume first that 1. holds. Let \(Y\) be a subspace of \(X\)
and consider the normalized weakly null sequences 
\((x_n)\) and \((x_n^*)\) in \(Y\) and \(Y^*\) respectively,
with \((x_n^*)\) \(\xi\)-convergent. Suppose that for some
\(\epsilon > 0\) and \(M =(m_i) \in [\mathbb{N}]\),
it was the case that \(x_{m_i}^*(x_{m_i}) > \epsilon\),
for all \(i \in \mathbb{N}\).
It follows that condition 2. of Theorem \ref{T:39}
is satisfied and therefore \((x_{m_i}^*)\)
admits a subsequence which is an \(\ell_1^{\xi}\)
spreading model in \(Y^*\). This contradicts
with Theorem \ref{T:3}, as \((x_n^*)\) is  \(\xi\)-convergent.
Hence \(Y\) satisfies the \(\xi\)-DP and 2. holds.

Now suppose that 2. holds. Let \((x_n)\) be a normalized
weakly null sequence in \(X\) admitting no subsequence which is
a \(c_0^{\xi}\) spreading model. In particular, no subsequence
of \((x_n)\) is equivalent to the unit vector basis of \(c_0\),
and thus by Corollary \ref{T:33} there exists \(M =(m_i) \in
[\mathbb{N}]\) such that the sequence
\((x_{m_i}^* |X_M)\) is weakly null in \(X_M^*\).
Where \((x_n^*)\) denotes the sequence of the biorthogonal functionals
of \((x_n)\).
Our assumption further yields that condition 2. of Theorem
\ref{T:39} fails for the space \(X_M\) and the weakly null sequences
\((x_{m_i})\) and \((x_{m_i}^*)|X_M\) in \(X_M\) and \(X_M^*\)
respectively. Thus condition 1. fails as well and so there
exists \(P=(p_i) \in [M]\)
so that \(x_{p_i}^*|X_P\) is \(\xi\)-convergent,
according to Theorem \ref{T:3}.
But \(X_P\) is \(\xi\)-DP and thus,
\(1= \lim_i x_{p_i}^*(x_{p_i}) =0\) which is absurd.
\end{proof}

\begin{remark}
Corollary \ref{T:40} can be viewed as a partial generalization of
Cembranos theorem, \cite{ce}, \cite{o3}, that a Banach space 
\(X\) is hereditarily Dunford-Pettis if, and only if,
every normalized weakly null sequence in \(X\) admits
a subsequence equivalent to the unit vector basis
of \(c_0\).
\end{remark}

\section{boundedly convexly complete sequences} \label{S:4}

This section is devoted to the proof of Theorem \ref{T:2}
which immediately yields Theorem \ref{T:1}. Our interest is
concentrated in weakly null sequences without 
boundedly convexly complete subsequences. In the next series of
lemmas, we describe their structure. 
We remark here that for the Schreier spaces
\(X_{\xi}\), \(\xi < \omega\), described in Section 
\ref{S:1}, it can be shown that they contain no 
boundedly convexly complete sequences.
However, an example of a weakly null
sequence \((y_n)\) in \(X_{\omega}\), is given in \cite{ao},
such that no convex block subsequence of \((y_n)\)
satisfies the weak Banach-Saks property. It turns out
that some subsequence of \((y_n)\), is boundedly 
convexly complete. We also mention that examples of
boundedly convexly complete sequences can be constructed
in \(C({\omega}^{\omega})\), the Banach space of
functions continuous on the ordinal interval
\([1,{\omega}^{\omega}]\) endowed by the order topology.

In the sequel, 
\(s=(x_n)\) will denote a normalized, bimonotone weakly null sequence
in the Banach space \(X\). We shall assume, without loss of
generality, that \(s\) satisfies the conclusion of Theorem
\ref{T:4}, for \(M=\mathbb{N}\) and \(\xi=0\).
That is, for every \(\delta > 0\), there exists \(C(\delta) > 0\)
such that if \((a_i)_{i=1}^n\) are scalars in \([-1,1]\), 
\(n \in \mathbb{N}\), and \(F \subset \{1, \cdots , n\}\)
with \(|a_i| \geq \delta\), for all \(i \in F\), then
\(\|\sum_{i \in F} a_i x_i \| \leq C(\delta) 
\|\sum_{i=1}^n  a_i x_i \|\).

\begin{notation}
\begin{enumerate}
\item Let \(t=(y_i)\) be a sequence in a Banach space and 
\(\boldsymbol{a} = (a_i)\) be a scalar sequence.
We let \(\boldsymbol{a} \cdot t\) denote the sequence
\((a_i y_i)\).
\item If \(M \in [\mathbb{N}]\), \(M=(m_i)\),
we let \(t_{M}\) denote the sequence \((y_{m_i})\).
\item Let \(P=(p_i)\) and \(Q=(q_i)\) be infinite subsets
of \(\mathbb{N}\) with \(P\) almost contained in \(Q\).
Then \({\boldsymbol{a}}_{Q|P}=(c_i)\) is the scalar sequence
obtained in the following manner:
Set \(c_i=0\), if \(q_i \notin P\). Then set
\(c_i=a_j\), if \(q_i = p_j\), for some \(j \in \mathbb{N}\).
\end{enumerate}
\end{notation}

\begin{Lem} \label{T:41}
Assume that \(s=(x_n)\) has no subsequence which is b.c.c.
There exist \(M_0 \in [\mathbb{N}]\) and 
\({\delta}_0 > 0\) such that for every \( L \in [M_0]\)
there exist a sequence of scalars \((a_m)_{m \in L}\) 
with \((a_m x_m)_{m \in L}\) strongly bounded by \(1\),
and a sequence \((F_i)_{i \in \mathbb{N}}\)
of consecutive subsets of \(L\) so that the 
following are satisfied for every \(i \in \mathbb{N}\):
\[ a_m \geq 0, \, \text{ for all }
   m \in F_i, \, \sum_{m \in F_i}a_m \leq 1, \,
\text{ and } 
\biggl \|\sum_{m \in F_i} a_m x_m \biggr \| > {\delta}_0.\]
\end{Lem}
\begin{proof}
We first observe that if
\(P=(p_i)\) and \(Q=(q_i)\) are infinite subsets
of \(\mathbb{N}\) with \(P\) almost contained in \(Q\), 
and \(\boldsymbol{a}\) is a scalar sequence,
then
\[ b({\boldsymbol{a}}_{Q|P} \cdot s_Q) \leq b(\boldsymbol{a} \cdot s_P), \,
\text{ and }
   c(\boldsymbol{a}, s_P) \leq c({\boldsymbol{a}}_{Q|P}, s_Q).\]
To prove the lemma, it suffices to find 
\(M_0 \in [\mathbb{N}]\), \({\delta}_1 > 0\) and 
\(1 \leq K < \infty\), so that for every 
\(L \in [M_0]\) there exists a scalar sequence
\(\boldsymbol{a}\) with
\(b(\boldsymbol{a} \cdot s_L) \leq K\) and
\(c(\boldsymbol{a}, s_L) > {\delta}_1 \).
Once this is accomplished, then
\(M_0\) and \({\delta}_0 = \frac{{\delta}_1}{2K}\),
satisfy the conclusion of the lemma.

We now claim that there exists \(N \in [\mathbb{N}]\)
and \( {\delta}_1 > 0\) such that
for all \(L \in [N]\) there exists a scalar sequence
\(\boldsymbol{a}\) with
\[b(\boldsymbol{a} \cdot s_L) < \infty,\, \text{ and }
  c(\boldsymbol{a}, s_L) > {\delta}_1 .\]
If this is not the case we construct a decreasing sequence
\((M_i)\) consisting of infinite subsets of \(\mathbb{N}\)
so that \(c(\boldsymbol{a}, s_{M_i}) \leq \frac{1}{i}\),
for every scalar sequence \(\boldsymbol{a}\) with
\(b(\boldsymbol{a} \cdot s_{M_i}) < \infty\).
Let now \(M \in [\mathbb{N}]\) be
almost contained in \(M_i\), for all \( i \in \mathbb{N}\),
and choose a scalar sequence \(\boldsymbol{a}\) such that
\[b(\boldsymbol{a} \cdot s_M) < \infty,\, \text{ and }
  c(\boldsymbol{a}, s_M) > 0.\]
It follows that
\[ b({\boldsymbol{a}}_{M_i|M} \cdot s_{M_i}) \leq 
b(\boldsymbol{a} \cdot s_M) < \infty ,\]
for all \(i \in \mathbb{N}\) and thus,
\[c(\boldsymbol{a}, s_M) \leq c({\boldsymbol{a}}_{M_i|M}, s_{M_i})
 \leq \frac{1}{i},\]
for all \(i \in \mathbb{N}\), which is absurd. Therefore our claim
holds.

We next claim that there exists \(M_0 \in [N]\) and 
\(K < \infty\) so that for every \(L \in [M_0]\), 
there exists a scalar sequence \(\boldsymbol{a}\) such that
\(b(\boldsymbol{a} \cdot s_L) \leq K \) and 
\(c(\boldsymbol{a}, s_L) > {\delta}_1\).
Again, were this claim false, we could choose a
decreasing sequence \((N_i)\) of infinite subsets of
\(N\) such that for all \(i \in \mathbb{N}\), if
\(c(\boldsymbol{a}, s_{N_i}) > {\delta}_1\),
for some scalar sequence \(\boldsymbol{a}\), then
\(b(\boldsymbol{a} \cdot s_{N_i}) > i\).
Now let \(T \in [N]\) be almost contained in each \(N_i\)
and choose a scalar sequence \(\boldsymbol{a}\)
such that 
\[ b(\boldsymbol{a} \cdot s_T) < \infty, \, \text{ and }\,
   c(\boldsymbol{a}, s_T) > {\delta}_1 .\]
It follows that
\[ c({\boldsymbol{a}}_{N_i|T}, s_{N_i}) \geq 
   c(\boldsymbol{a}, s_T) > {\delta}_1,\]
for all \(i \in \mathbb{N}\), and hence
\[ b(\boldsymbol{a} \cdot s_T) \geq 
b({\boldsymbol{a}}_{N_i|T} \cdot s_{N_i}) > i ,\]
for all \(i \in \mathbb{N}\) which is absurd.
\end{proof}

\begin{remark}
If \(L \in [M_0]\), \((F_i)_{i=1}^{\infty}\) and \((a_m)_{m \in L}\),
are as in the conclusion of Lemma \ref{T:41},
then we shall call the sequence
\((\sum_{m \in F_i}a_m x_m)_{i \in \mathbb{N}}\),
a sub-convex block subsequence of \((x_i)\)
supported by \(L\) and satisfying the conclusion of
Lemma \ref{T:41}.
\end{remark}

\begin{Lem} \label{T:42}
Let \(s=(x_n)\) be a weakly null sequence having no
subsequence which is b.c.c. 
Let \(M_0 \in [\mathbb{N}]\) and \({\delta}_0 > 0\)
be as in the conclusion of Lemma \ref{T:41} applied on \(s\).
Suppose \(\alpha < {\omega}_1\)
is a limit ordinal and denote by \(({\alpha}_n + 1)\) 
the sequence of ordinals associated to \(\alpha\).
Assume that for every \(n \in \mathbb{N}\) and every
\(N \in [M_0]\) there
exists \(M \in [N]\) such that \((x_m)_{m \in M}\)
is an \(l_{1}^{{\alpha}_n + 1}\) spreading model.
Then, for every \(N \in [M_0]\), there 
exists \(M \in [N]\) such that \((x_m)_{m \in M}\)
is an \(l_{1}^{\alpha}\) spreading model.
\end{Lem}
\begin{proof}
Let \(\mathcal{A}\) denote the set of those sub-convex
combinations of the sequence \(s\) which are members of
a (not necesserilly the same)
sub-convex block subsequence of \(s\) that satisfies the
conclusion of Lemma \ref{T:41}.

We shall apply Corollary \ref{T:212} to the sequence \(s\)
and the family \(\mathcal{A}\) in order to obtain
\(T \in [N]\) so that \(S_{\alpha}(T) \subset 
{\mathfrak{F}}_{\frac{{\delta}_0}{2}}\).
Corollary \ref{T:36} will then yield that for some
\(M \in [T]\), \((x_m)_{m \in M}\) is an 
\(l_{1}^{\alpha}\) spreading model.
To this end, it suffices to show that for every \(i_0 \in 
\mathbb{N}\), every
\(\epsilon > 0\) and \(L \in [N]\), there exists \(x \in \mathcal{A}\)
supported by \(L\) and such that \(\|x\|_{{\alpha}_{i_0}} 
< \epsilon\). Suppose this is not the case and choose
according to the hypothesis
\(P \in [L]\) so that \((x_m)_{m \in P}\) is an
\(l_{1}^{{\alpha}_{i_0} + 1}\) spreading model with constant \(K\).
Without loss of generality, by Theorem \ref{T:6}, we can
assume that \((x_m)_{m \in P}\) is 
\(S_{{\alpha}_{i_0} + 1}\)
unconditional with constant \(C\).
Next, choose by Lemma \ref{T:24}, \(Q \in [P]\)
such that if \(F \in S_{{\alpha}_{i_0} + 1}[Q]\), 
then \(F \setminus \{\min F\} \in S_{{\alpha}_{i_0} + 1}(P)\).

Since \(s\) has no b.c.c. subsequences, there exist, by Lemma
\ref{T:41}, a sequence of scalars
\((a_q)_{q \in Q}\) and a sequence \((F_i)_{i=1}^{\infty}\)
of consecutive subsets of \(Q\) with
\((a_q x_q)_{q \in Q}\) strongly bounded by \(1\),  
such that for all \( i \in \mathbb{N}\),
\[ a_q \geq 0, \text{ for all } \, q \in F_i, \,
   \sum_{q \in F_i}a_q \leq 1, \, \text{ and }
\|\sum_{q \in F_i}a_q x_q \| > {\delta}_0 .\]
However, \(\sum_{q \in F_i}a_q x_q \in \mathcal{A}\),
for all \( i \in \mathbb{N}\), and moreover they are
supported by \(L\).
Thus, there exists, for all \( i \in \mathbb{N}\),
\[G_i \subset F_i, \, 
 G_i \in S_{{\alpha}_{i_0}} \text{ so that } \,
\sum_{q \in G_i} a_q \geq \epsilon .\]
Next choose \(p \in \mathbb{N}\) such that 
\((p-1)K \epsilon > C \). Choose also
\(j \in \mathbb{N}\) with \( p < \min F_j\). Then
\[ \bigcup_{l=j + 1}^{j + p}G_l \in 
S_{{\alpha}_{i_0} + 1}[Q] \, \text{ and thus } \,
   \bigcup_{l=j + 2}^{j + p}G_l \in S_{{\alpha}_{i_0} + 1}(P).\]
Therefore,
\[
\biggl \| \sum_{ q \in \bigcup_{l=j + 2}^{j + p}G_l}
            a_q x_q \biggr \| \geq
 K \sum_{ q \in \bigcup_{l=j + 2}^{j + p}G_l} a_q \geq
 K(p-1) \epsilon > C. \]
However,
\[\biggl \| \sum_{ q \in \bigcup_{l=j + 2}^{j + p}G_l}
            a_q x_q \biggr \| \leq 
   C \biggl \| \sum_{q \leq \max G_{j +p}} a_q x_q \biggr \|
   \leq C.\]
This contradiction completes the proof of the lemma.
\end{proof}

\begin{Lem} \label{T:43}
Let \(s=(x_n)\) be a weakly null sequence having no subsequence
which is b.c.c. There exist \(M \in [\mathbb{N}]\),
\(\xi < {\omega}_1\) and \(\delta > 0\) so that the following
are satisfied:
\begin{enumerate}
\item \((x_m)_{m \in M}\) is an \(\ell_1^{\xi}\) spreading model
yet no subsequence of \((x_m)_{m \in M}\) is an
\(\ell_1^{\xi + 1}\) spreading model.
\item For every \(N \in [M]\), there exist \(L \in [N]\)
and a sequence of scalars \((a_m)_{m \in N}\) with
\((a_m x_m)_{m \in N}\) strongly bounded by \(1\)
so that for all \(i \in \mathbb{N}\)
\[ a_m \geq 0, \, \text{ for all }
   m \in F_i^{\xi}(L), \, 
\sum_{m \in F_i^{\xi}(L)} a_m  \leq 1, \,
\text{ and } 
\biggl \|\sum_{m \in F_i^{\xi}(L)}a_m x_m \biggr \| > 
 \delta .\]
\end{enumerate}
\end{Lem}
\begin{proof}
Let \(M_0 \in [\mathbb{N}]\) and \({\delta}_0 > 0\)
be as in the conclusion of Lemma \ref{T:41}
applied on \(s\).
\begin{align}
\text{ Let } \zeta = &\min \{ \alpha <{\omega}_1 : \,
\exists \, P \in [M_0] \, \text{ such that } \, \forall
L \in [P], \, \notag \\
&(x_n)_{n \in L} \text{ is not an } \ell_1^{\alpha}
\text{ spreading model } \}.\notag
\end{align}
Lemma \ref{T:42} yields that \(\zeta\) is not a limit ordinal and thus
\(\zeta = \xi + 1\), for some countable ordinal \(\xi\).
Choose now \(P \in [M_0]\) so that no subsequence of
\((x_n)_{ n \in P}\) is an \(\ell_1^{\xi + 1}\) spreading model.
Since \( \xi < \zeta\), there exists \(Q \in [P]\) so that
\((x_n)_{ n \in Q}\) is an \(\ell_1^{\xi}\) spreading model
and of course no subsequence of \((x_n)_{ n \in Q}\)
is an \(\ell_1^{\xi + 1}\) spreading model.

It follows by Corollary
\ref{T:36}, that for every \(L \in [Q]\),
the family \(S_{\xi + 1}(L)\) is not contained in
\({\mathfrak{F}}_{\frac{{\delta}_0}{2}}\).
Using Lemma \ref{T:29} as we did in the proof of Theorem
\ref{T:211}, we find \(M \in [Q]\) and 
\(n_0 \in \mathbb{N}\), so that
\[ \bigcup_{i=1}^{n_0} F_{i}^{\xi}(L) \notin 
  {\mathfrak{F}}_{\frac{{\delta}_0}{2}}\]
for all \(L \in [M]\).
It follows now, by Lemma \ref{T:26},
that every member of the family
\({\mathfrak{F}}_{\frac{{\delta}_0}{2}}[M]\)
can be expressed as the union of at most \(n_0\)
consecutive \(S_{\xi}\) sets.

We are going to show that \(M\) and
\( \delta = \frac{{\delta}_0}{2n_0}\) satisfy 2.
Indeed, let \(N \in [M]\) and choose
\[ n_1 < A_1 < n_2 < A_2 < \cdots \]
so that for all \( i \in \mathbb{N}\), 
\(n_i \in N\) and \(A_i\) is a maximal \(S_{\xi}\)
subset of \(N\).
Set \(N_1= \{n_i : i \in \mathbb{N}\}\).
By Lemma \ref{T:41}, there exist a sequence of scalars
\((a_m)_{m \in N_1}\) and a sequence of consecutive
subsets of \(N_1\), \((F_i)_{i=1}^{\infty}\) with
\((a_m x_m)_{m \in N_1}\) strongly bounded by \(1\),
so that for every \(i \in \mathbb{N}\)
\[ a_m \geq 0, \, \text{ for all }
   m \in F_i, \, \sum_{m \in F_i}a_m \leq 1, \,
\text{ and } 
\biggl \|\sum_{m \in F_i} a_m x_m \biggr \| > {\delta}_0.\]
We next choose, for every \(i \in \mathbb{N}\),
a functional \(x_i^* \in B_{X^*}\) such that 
\[ \biggl | \sum_{m \in F_i} a_m x_i^*(x_m) \biggr | > {\delta}_0 , \]
and let
\[ G_i = \{ m \in F_i : |x_i^*(x_m)| \geq \frac{{\delta}_0}{2}\}.\]
It follows now, as \(G_i\) belongs to 
\({\mathfrak{F}}_{\frac{{\delta}_0}{2}}[M]\), 
that there exists, for every \(i \in \mathbb{N}\),
\[H_i \subset G_i, \, H_i \in S_{\xi},\,
\text{ such that }\,
\biggl | \sum_{m \in H_i} a_m x_i^*(x_m) \biggr | >
\frac{{\delta}_0}{2 n_0} = \delta . \]
By extending \(H_i\) to a maximal \(S_{\xi}\) subset of
\(H_i \bigcup A_{\max F_i}\), if necessary, we can assume
without loss of generality, that \(H_i\) itself is a 
maximal \(S_{\xi}\) subset of \(N\), for all \(i \in \mathbb{N}\).
We also extend the sequence \((a_m)_{m \in N_1}\) to
a scalar sequence \((a_m)_{m \in N}\) in the
obvious manner.

Concluding, there exist \(H_1 < H_2 < \cdots \) maximal 
\(S_{\xi}\) subsets of \(N\) and a scalar sequence
\((a_m)_{m \in N}\) with 
\((a_m x_m)_{m \in N}\) strongly bounded by \(1\),
so that for every \(i \in \mathbb{N}\)
\[ a_m \geq 0, \, \text{ for all }
   m \in H_i, \, \sum_{m \in H_i}a_m \leq 1, \,
\text{ and } 
\biggl \|\sum_{m \in H_i} a_m x_m \biggr \| > \delta .\]
Finally, choose \(L \in [N]\) such that 
\(F_i^{\xi}(L)=H_i\), for all \(i \in \mathbb{N}\),
and we are done.
\end{proof} 

\begin{Def}
Let \(M \in [\mathbb{N}]\), \(L \in [M]\), 
\(n \in \mathbb{N}\) and \(\delta > 0\).
The scalar sequence \((a_m)_{m \in M}\)
is called \(\xi\)-faithful 
for \((\delta, L , n )\), provided the following 
conditions hold:
\begin{enumerate}
\item \(|a_m| \leq \frac{\delta}{2}\), for all \(m \in M\).
\item \((a_m x_m)_{m \in M}\) is strongly bounded by \(1\).
\item For every \(i \leq n\) we have that 
\[ a_m \geq 0, \, \text{ for all }
   m \in F_i^{\xi}(L), \, 
\sum_{m \in F_i^{\xi}(L)} a_m  \leq 1, \,
\text{ and } 
\biggl \|\sum_{m \in F_i^{\xi}(L)}a_m x_m \biggr \| \geq 
 \delta .\]
\end{enumerate}
\end{Def}

\begin{Lem} \label{T:44}
Let \(s=(x_n)\) have no subsequence which is b.c.c.
and satisfying the conclusion of Lemma \ref{T:43}
for some \(M \in [\mathbb{N}]\), \(\xi < {\omega}_1\)
and \(\delta > 0\).
\begin{enumerate}
\item If \(\xi =0\), then some subsequence of
\((x_m)_{m \in M}\) is equivalent to the unit vector basis
of \(c_0\).
\item If \(\xi \geq 1\), then for every \(N \in [M]\)
there exists \(P \in [N]\) so that for all 
\(L \in [P]\) and \(n \in \mathbb{N}\)
there exists a scalar sequence \((a_m)_{m \in N}\)
which is \(\xi\)-faithful for \((\delta, L , n )\).
\end{enumerate}
\end{Lem}
\begin{proof}
If \(\xi =0\), then according to Lemma \ref{T:43}
there exist a scalar sequence \((a_m)_{m \in M}\)
and \(L \in [M]\), \(L=(l_i)\), so that
\((a_m x_m)_{m \in M}\) is strongly bounded by \(1\)
and \(a_m > \delta\), for all \(m \in L\).
It follows now that if \(x^* \in B_{X^*}\), then for all
\(k \in \mathbb{N}\), 
\[ \sum_{i=1}^k |x^*(x_{l_i})| \leq \frac{2C(\delta)}{\delta}\]
and thus \((x_{l_i})\) is equivalent to the unit vector basis of
\(c_0\).

Assume now that \(\xi \geq 1\) and let \(N \in [M]\).
It follows that no subsequence of
\((x_m)_{m \in M}\) is equivalent to the unit vector basis
of \(c_0\) and therefore, by Corollary \ref{T:33}, we
obtain \(N_0 \in [N]\) such that \((x_m)_{m \in N_0}\)
is semi-boundedly complete.

We next choose, according to Lemma \ref{T:43}, a scalar
sequence \((a_m)_{m \in N_0}\) and \(L_0 \in [N_0]\)
with \((a_m x_m)_{m \in N_0}\) strongly bounded
by \(1\) so that for all \(i \in \mathbb{N}\),
\[ a_m \geq 0, \, \text{ for all }
   m \in F_i^{\xi}(L_0), \, 
\sum_{m \in F_i^{\xi}(L_0)} a_m  \leq 1, \,
\text{ and } 
\biggl \|\sum_{m \in F_i^{\xi}(L_0)}a_m x_m \biggr \| \geq 
 \delta .\]
Since \(\lim_{m \in N_0}a_m =0\), there exists \(i_0 \in \mathbb{N}\)
such that \(|a_m| < \frac{\delta}{2}\), for all \( m \in N_0\),
\( m > \max F_{i_0}^{\xi}(L_0)\). Set
\[ L = \bigcup_{i > i_0}F_{i}^{\xi}(L_0) \]
and extend the sequence \((a_m)_{m \in N_0}\) 
to a sequence \((a_m)_{m \in N}\) in the obvious manner.
Evidently,
\((a_m)_{m \in N}\) is \(\xi\)-faithful for 
\((\delta, L , n)\), for all \(n \in \mathbb{N}\).
We next consider the set
\[\mathcal{D} = \{ L \in [N]:\, \forall \, n \in \mathbb{N} \quad
\exists \, (a_m)_{m \in N} \quad \xi \text{ -faithful for }
(\delta, L , n)\}\]
which is of course closed in \([N]\) and therefore Ramsey.
Our previous argument yields the existence of \(P \in [N]\)
such that \([P] \subset \mathcal{D}\). Clearly \(P\) satisfies
the conclusion of part 2. of this lemma. 
\end{proof}
\begin{Lem} \label{T:45}
Let \(s = (x_n)\) be a weakly null sequence and assume that
\((x_n^*)\), the sequence of functionals biorthogonal to
\((x_n)\) is weakly null in \([(x_n)]^*\).
Let \(\xi < {\omega}_1\), \(M \in [\mathbb{N}]\),
\(\epsilon > 0\) and \(\delta > 0\).
There exists \(N \in [M]\) satisfying the following property:
If \(L \in [N]\), \(n \in \mathbb{N}\)
and there exists a scalar sequence \((a_m)_{m \in M}\)
which is \(\xi\)-faithful for \((\delta, L , n)\),
then there exists \((b_m)_{m \in M}\) which is
\(\xi\)-faithful for \((\delta, L , n)\)
and so that
\[ \sum_{ m \in N \setminus \bigcup_{i=1}^{n}
   F_{i}^{\xi}(L)} |b_m| < \epsilon . \]
\end{Lem}
\begin{proof}
The proof is similar to those of Lemmas \ref{T:32} and 
\ref{T:34}. 
Choose first a sequence \(({\epsilon}_i)_{i=0}^{\infty}\)
of positive scalars such that 
\(\sum_{i=0}^{\infty} {\epsilon}_i < \epsilon\).
Using the same notation and terminology as in Lemma \ref{T:32},
let \(n \in \mathbb{N} \bigcup \{0\}\)
and \(F \subset T_n\).
The \(n\)-tuple of elements of \(M\)
\((r_i)_{i \in T_n}\) and the infinite
subset \(L\) of \(M\), 
\(L=(l_i)\), 
are said to satisfy property \((F-A_n)\)
provided that \(r_n < l_1\), if \(n \geq 1\),
and the following statement holds:

If \( k \in \mathbb{N}\) and there exists a scalar sequence
\((a_m)_{m \in M}\) which is \(\xi\)-faithful for
\((\delta , 
\{r_i : i \in T_n \setminus F\} \bigcup \{l_i : i \geq 2\},
k)\),
then there exists \((b_m)_{m \in M}\) which is \(\xi\)-faithful for
\((\delta , 
\{r_i : i \in T_n \setminus F\} \bigcup \{l_i : i \geq 2\},
k)\) and such that
\[\sum_{i \in F}|b_{r_i}| + |b_{l_1}| < 
\sum_{i=0}^{n} {\epsilon}_i .\]
Let us also say that \((r_i)_{i \in T_n}\) and \(L\)
satisfy property \((A_n)\), if they satisfy property
\((F-A_n)\), for all \(F \subset T_n\).

We shall inductively construct an increasing sequence
\((r_n)_{n=1}^{\infty}\) of elements of \(M\)
and a decreasing sequence \((M_n)_{n=0}^{\infty}\)
of infinite subsets of \(M\) with
\(r_n \in M_{n-1}\), if \(n \geq 1\), so that if 
\(n \in \mathbb{N} \bigcup \{0\}\) and \(L \in [M_n]\), then
\((r_i)_{i \in T_n}\) and \(L\) satisfy property \((A_n)\).

The first inductive step is similar to the general one 
and so we shall not explicit it. Now we assume that
\(r_1 < \cdots < r_n\) and \(M_0 \supset \cdots \supset M_n\)
have been constructed with \(r_i \in M_{i-1}\), if \(i \geq 1\),
so that if \(i \leq n\), and \(L \in [M_i]\) then
\((r_j)_{j \in T_i}\) and \(L\) satisfy property \((A_i)\).
Let \(r_{n+1} = \min M_n\) and fix \(F \subset T_{n+1}\).
We define
\[ {\mathcal{\Delta}}_{F}= \{ L \in [M_n]:\,  (r_i)_{i \in T_{n+1}}
\, \text{ and } L \, \text{ satisfy } \, (F-A_{n+1})\}. \] 
Clearly, \({\mathcal{\Delta}}_{F}\) is closed in \([M_n]\)
and therefore Ramsey. Suppose that for some \(P \in [M_n]\),
\(P=(p_i)\), we had that \([P] \bigcap {\mathcal{\Delta}}_{F}=
\emptyset\). Let \(q \in \mathbb{N}\) and set
\[ L_j = \{p_j\} \bigcup \{p_i : \, i > q\}, \,
\text{ for all } j \leq q . \]
Since \(L_j \notin {\mathcal{\Delta}}_{F}\), for all 
\(j \leq q\), there exist integers \((k_j)_{j=1}^{q}\)
as well as scalar sequences \((a_{m}^j)_{m \in M}\),
\(j \leq q\), so that letting
\(R= \{ r_i : \, i \in T_{n+1} \setminus F\}
\bigcup \{p_i : \, i > q\}\) we have that for all
\(j \leq q\),
\[ (a_{m}^{j})_{m \in M} \, \text{ is } \, \xi 
\text{ -faithful for } \, 
(\delta, R , k_j), \]
and moreover, if \((b_m)_{m \in M}\) is a 
\(\xi\)-faithful scalar sequence for 
\((\delta, R , k_j)\) then
\[ \sum_{i \in F} |b_{r_i}| + |b_{p_j}| \geq
\sum_{i=0}^{n +1} {\epsilon}_i .\]
Let \(k_0 = \max \{ k_j : \, j \leq q\}\)
and observe that any scalar sequence 
\((b_m)_{m \in M}\) which is \(\xi\)-faithful
for \((\delta, R , k_0)\), is also
\(\xi\)-faithful for \((\delta, R , k_j)\),
for all \(j \leq q\).
Next, let \(t = \max F\) and note that
\[R= \{r_i :\,  i \in T_{t-1} \setminus F\} \bigcup
\{ r_i : \, t < i \leq n+1\} \bigcup
\{p_i : \, i > q\}.\]
By the induction hypothesis, since \((r_i)_{i=1}^{t-1}\)
and \(\{ r_i : \, t \leq i \leq n+1\} \bigcup
\{p_i : \, i > q\}\)
satisfy property \((A_{t-1})\),
there exists a scalar sequence 
\((b_m)_{m \in M}\) which is \(\xi\)-faithful
for \((\delta, R , k_0)\) and such that
\[\sum_{i \in F \setminus \{t\}} |b_{r_i}| + 
|b_{r_t}| < \sum_{i=0}^{t-1} {\epsilon}_i .\]
Thus,
\[\sum_{ i \in F} |b_{r_i}| < 
\sum_{i=0}^{n} {\epsilon}_i .\]
It follows now by our previous observation, that
\((b_m)_{m \in M}\) is \(\xi\)-faithful for
\((\delta, R , k_j)\),
for all \(j \leq q\) and thus,
\[ \sum_{i \in F} |b_{r_i}| + |b_{p_j}| \geq
\sum_{i=0}^{n +1} {\epsilon}_i , \, 
\text{ for all } j \leq q .\]
Hence, \(|b_{p_j}| \geq {\epsilon}_{n+1}\),
for all \( j \leq q\).

Now choose \(m_0 > p_q\), \( m_0 \in M\).
Then, \(|x_{p_j}^{*}( \sum_{ m \leq m_0} b_m x_m )|
\geq {\epsilon}_{n+1}\), for all
\(j \leq q\). 
But \(\|\sum_{ m \leq m_0} b_m x_m \| \leq 1\)
and \(q\) is arbitrary and so \((x_{p_i}^*)\)
is not weakly null in \([(x_i)]^*\),
contradicting our assumption.

Concluding, there exists \(L \in [M_n]\)
such that \([L] \subset {\mathcal{\Delta}}_{F}\).
By repeating the previous argument successively over
all possible subsets of \(T_{n+1}\),
we obtain \(M_{n+1} \in [M_n]\) such that
\([M_{n+1}] \subset {\mathcal{\Delta}}_{F}\),
for all \(F \subset T_{n+1}\).
The inductive construction is now complete.
Set \(N =(r_i)\). Let \(L \in [N]\)
and \(k \in \mathbb{N}\).
Let also \((a_m)_{m \in M}\) be a scalar sequence
which is \(\xi\)-faithful for \((\delta , L , k)\).
Suppose that \(r_p = \max F_{k}^{\xi}(L)\)
and let \(G_q= \{ i \in T_q : \, r_i \notin 
\bigcup_{j=1}^{k} F_{j}^{\xi}(L)\}\), for all \(q \geq p\).
Our construction yields that \((r_i)_{i=1}^{q}\)
and \(\{r_i : \, i > q\}\), satisfy
\((G_q - A_q)\), for all \(q \geq p \).
We also have, by stability, that for all \(q \geq p\),
\[ F_{j}^{\xi}(L) = F_{j}^{\xi}(R_q), \, \text{ for all }
j \leq k \]
where, \(R_q = \{r_i : i \in T_q \setminus G_q\}
\bigcup \{r_i : \, i > q+1\}\),
and therefore there exists a scalar sequence 
\((b_m^q)_{m \in M}\) which is 
\(\xi\)-faithful for \((\delta , L , k)\)
and such that
\[\sum_{i \in G_q}|b_{r_i}^q| < \sum_{i=0}^{q} {\epsilon}_i < \epsilon .\]
Thus, 
\[\sum_{ m \in N \setminus \bigcup_{i=1}^{k} F_{i}^{\xi}(L) ,
\, m \leq r_q} |b_{m}^q| < \epsilon.\]
Finally let \((b_m)_{m \in M}\) be any cluster point of
the sequence \(\bigl( (b_m^q)_{m \in M} \bigr )_{q \geq p}\)
in \([-1,1]^M\). Evidently, this is the desired sequence.
\end{proof}
\begin{Lem} \label{T:46}
Let \(t= (y_i)\) be a sequence in \(B_X\), \(\xi < {\omega}_1\)
and \(\mu\) a finitely supported probability measure on
\(\mathbb{N}\).
Assume there exist \(x^* \in B_{X^*}\) and \(\epsilon > 0 \)
so that \(x^*(\mu \cdot t) \geq 2 \epsilon\).
Assume further that there exists \(l \in \mathbb{N}\),
\(l < supp \mu\), such that
\(l \|\mu\|_{ {\zeta}_i} < \epsilon\), for every
\(i \leq l\),
where \(({\zeta}_i + 1)\) is the sequence of ordinals associated
to \(\xi\). Let
\[E = \{ n \in supp \mu : x^* (y_n) \geq \frac{1}{2} 
x^*(\mu \cdot t) \}. \]
Then, there exists a maximal \(S_{\xi}\) set containing \(l\)
and contained in \(\{l\} \bigcup E\).
\end{Lem}
\begin{proof}
Assume on the contrary, that no subset of \(\{l\} \bigcup E\)
containing \(l\) is a maximal \(S_{\xi}\) set.
We then claim that \(\{l\} \bigcup E\) belongs to \(S_{\xi}\).
Indeed, suppose that \(\{l\} \bigcup E
= \{m_1, \cdots, m_k \}\), \(m_1 =l\), and choose \(r \leq k\)
maximal with respect to \(\{m_1, \cdots, m_k \} \in S_{\xi}\).
If \(r < k \), then Lemma \ref{T:25} yields that
\(\{m_1, \cdots, m_r \}\) is a maximal \(S_{\xi}\) set
contradicting our assumption. Thus \(r=k\) as claimed.

There exists now \( i \leq l\) such that \(\{l\} \bigcup E\)
belongs to \(S_{{\zeta}_i + 1}\). It follows that we can find
\(p \leq l\) and consecutive \(S_{{\zeta}_i}\) sets
\((A_j)_{j=1}^{p}\) so that
\[ \{l\} \bigcup E = \bigcup_{j=1}^p A_j .\]
But now,
\begin{align}
x^*(\mu \cdot t) = &\sum_{n \in E} \mu(n) x^*(y_n) +
\sum_{n \in supp \mu \setminus E} \mu(n) x^*(y_n) \notag \\
&\leq \sum_{j=1}^p \mu(A_j) + \frac{1}{2}x^*(\mu \cdot t) \notag \\
&\leq p\|\mu\|_{{\zeta}_i}  + \frac{1}{2}x^*(\mu \cdot t) \notag 
\end{align}
and thus \(x^*(\mu \cdot t) < 2 \epsilon\). This contradiction
completes the proof of the lemma.
\end{proof}

We are now ready for the proof of the main result of this paper.
\begin{proof}[Proof of Theorem \ref{T:2}]
Assume that \(s\) has no subsequence which is b.c.c. Choose
\(M \in [\mathbb{N}]\), \(\xi < {\omega}_1\) 
and \(\delta > 0\) satisfying the conclusion of Lemma \ref{T:44}
applied on \(s\).
If \(\xi =0\) we are done since some subsequence of \(s\) is
equivalent to the unit vector basis of \(c_0\).

Assume now that \(\xi \geq 1\). Choose according to Corollary
\ref{T:33}, \(M_1 \in [M]\) so that the sequence 
\((x_m^*)_{m \in M_1}\) of functionals biorthogonal to
\((x_m)_{m \in M_1}\) is weakly null in \([(x_m)_{m \in M_1}]^*\). 
Next choose, according to Lemma \ref{T:24}, \(M_2 \in [M_1]\)
such that
\begin{equation} \label{E:2}
F \setminus \{\min F\} \in S_{\xi}(M_1),\, 
\text{ for all } \, F \in S_{\xi}[M_2] . 
\end{equation}
Let \(0 < \epsilon < 1\) and choose a sequence of positive
scalars \(({\epsilon}_n)\) such that 
\(\sum_{n=1}^{\infty}{\epsilon}_n < \epsilon\). Choose also
\(\lambda \geq \frac{4}{\delta} (1 + 4 \epsilon) \).
Lemma \ref{T:32} now yields \(M_3 \in [M_2]\) such that for
every \(L \in [M_3]\) and \(n \in \mathbb{N}\), if there
exists \(x^* \in B_{X^*}\) which is \(\xi\)-good for
\((\delta, L , n)\), then there exists 
\(y^* \in B_{X^*}\), \(\xi\)-good for \((\delta, L , n)\)
and such that
\begin{equation} \label{E:3}
\sum_{ m \in M_3 \setminus \bigcup_{i=1}^{n}
       supp {\xi}_{i}^{L} }
 |y^*(x_m)| < \epsilon . 
\end{equation}
We continue our choice of infinite subsets of \(M\) by choosing
\(M_4 \in [M_3]\) according to Lemma \ref{T:45}. Thus,
for every \(L \in [M_4]\) and \(n \in \mathbb{N}\), if there
exists a scalar sequence \((a_m)_{m \in M_3}\)
which is \(\xi\)-faithful for \((\delta, L , n)\), then there exists
\((b_m)_{m \in M_3}\), 
\(\xi\)-faithful for \((\delta, L , n)\), and such that
\begin{equation} \label{E:4}
\sum_{m \in M_4 \setminus \bigcup_{i=1}^{n}
      F_{i}^{\xi}(L)}
|b_m| < \epsilon .
\end{equation}
Finally, choose \(M_5 \in [M_4]\) according to Lemma \ref{T:44}
applied for ``\(N\)'' \(= M_4\). It follows 
now, by (\ref{E:4}), that for
every \(L \in [M_5]\) and \(n \in \mathbb{N}\), there exists
a scalar sequence \((a_m)_{m \in M_3}\) 
which is 
\(\xi\)-faithful for \((\delta, L , n)\)
and such that
\begin{equation} \label{E:5}
\sum_{m \in M_4 \setminus \bigcup_{i=1}^{n}
      F_{i}^{\xi}(L)}
|a_m| < \epsilon .
\end{equation}

Let now \(Q \in [M_5]\). Let also
\(({\zeta}_n + 1)\) be the sequence of ordinals associated to
\(\xi\).
Repeated applications of Proposition \ref{T:215}, 
now yield an increasing
sequence of elements of \(Q\), \((l_i)\), and a sequence
\((Q_i)\) of infinite subsets of \(Q\), so that
\[ l_1 < supp {\xi}_{1}^{Q_1} < l_2 < supp {\xi}_{1}^{Q_2} < 
\cdots \]
and 
\[
 l_i \|{\xi}_{1}^{Q_i}\|_{\alpha} < {\epsilon}_i , \, 
\text{ for all } \, \alpha \in 
\{ {\zeta}_m : m \leq l_i \} \bigcup \{0\} 
\]
for all \(i \in \mathbb{N}\).
We thus obtain, by stability, \(P \in [Q]\) such that
\({\xi}_{i}^{P}= {\xi}_{1}^{Q_i}\), for all \(i \in \mathbb{N}\),
and therefore,
\begin{align} 
&l_1 < supp {\xi}_{1}^{P} < l_2 < supp {\xi}_{2}^{P} < 
\cdots \notag \\
&l_i \|{\xi}_{i}^{P}\|_{\alpha} < {\epsilon}_i , \, 
\text{ for all } \, \alpha \in 
\{ {\zeta}_m : m \leq l_i \} \bigcup \{0\} \, \text{ and }
i \in \mathbb{N}. \label{E:6}
\end{align} 
Let \(K\) be the \(\ell_1^{\xi}\) spreading model constant of 
\((x_m)_{m \in M}\). By our choice of \(M_2 \in [M_1]\),
(\ref{E:2}) and (\ref{E:6}) yield that
\[\|{\xi}_{i}^{P} \cdot s\| \geq K - \frac{{\epsilon}_i}{l_i}
(K + 1)\]
and hence \(({\xi}_{i}^{P} \cdot s)\) is semi-normalized.

We now {\em claim \/} that
\[ \biggl \| \sum_{i \in F} {\xi}_{i}^{P} \cdot s
        \biggr \| \leq 2 {\lambda}, \]
for all \(F \in [\mathbb{N}]^{< \infty}\). Our claim of course
implies that \(({\xi}_{i}^{P} \cdot s)\) is equivalent to
the unit vector basis of \(c_0\).
Were our claim false, there would exist 
\(q \in \mathbb{N}\), integers 
\(i_1 < \cdots <  i_q \), and \(y^* \in B_{X^*}\)
such that
\[ y^*( {\xi}_{i_n}^{P} \cdot s ) \geq 0, \, \text{ for }
n \leq q, \, \text{ and } \, 
\sum_{n=1}^{q} y^*( {\xi}_{i_n}^{P} \cdot s ) > \lambda . \]
Since \(P \in [M_3]\), (\ref{E:3}) yields \(x^* \in B_{X^*}\)
such that
\begin{align} 
&x^*( {\xi}_{i_n}^{P} \cdot s ) \geq 0, \, \text{ for }
n \leq q, \, 
\sum_{n=1}^{q} x^*( {\xi}_{i_n}^{P} \cdot s ) > \lambda \text{ and } 
\notag \\
&\sum_{ m \in M_3 \setminus \bigcup_{n=1}^{q}
       supp {\xi}_{i_n}^{P} }
 |x^*(x_m)| < \epsilon . \label{E:7}
\end{align}
Next we set \(I = \{ n \leq q : \, x^*( {\xi}_{i_n}^{P} \cdot s )
 \geq 2 {\epsilon}_{i_n}\}\). It follows that
\begin{equation} \label{E:8}
\sum_{n \in I} x^*( {\xi}_{i_n}^{P} \cdot s ) > 
\lambda - 2 \epsilon .
\end{equation}
We now fix \(n \in I\) and let
\[ E_n = \{ m \in supp {\xi}_{i_n}^{P} : \, 
x^*(x_m) \geq \frac{1}{2} x^*( {\xi}_{i_n}^{P} \cdot s ) \}. \]
Because of (\ref{E:6}), Lemma \ref{T:46} 
yields the existence of a maximal
\(S_{\xi}\) subset \(D_n\) of \(\{l_{i_n}\} \bigcup E_n\)
such that \(\min D_n = l_{i_n}\).
Let \(m_0 = \max D_{\max I}\) and set
\[ L = \bigcup_{ n \in I} D_n \bigcup 
\{ m \in Q : \, m > m_0 \}.\] 
Since \(L \in [M_5]\), we obtain through (\ref{E:5}),
a scalar sequence \((a_m)_{m \in M_3}\) which is 
\(\xi\)-faithful for \((\delta , L , |I|)\) and such that
\begin{equation} \label{E:9}
\sum_{m \in M_4 \setminus \bigcup_{ n \in I} D_n}
|a_m| < \epsilon .
\end{equation}
We recall here that for every \(n \in I\),
\[ a_m \geq 0 , \, \text{ for all } \, m \in D_n , \,
   \sum_{ m \in D_n} a_m \leq 1 , \quad
\text{ and } \, 
\biggl \| \sum_{m \in D_n} a_m x_m \biggr \| > \delta .\] 
Of course, \((a_m x_m)_{m \in M_3}\) is strongly bounded by
\(1\) and \(|a_m| \leq \frac{\delta}{2}\), for all
\(m \in M_3\). Therefore,
\begin{equation} \label{E:10}
\sum_{m \in D_n \setminus \{l_{i_n}\}} 
a_m > \frac{\delta}{2}, \, \text{ for all } \,
n \in I .
\end{equation}
Our construction yields that
\begin{equation} \label{E:11}
\sum_{m \in M_3 \setminus \bigcup_{ n \in I} D_n}
|a_m||x^*(x_m)| < 2 \epsilon .
\end{equation}
Indeed, 
\begin{align}
\sum_{m \in M_3 \setminus \bigcup_{ n \in I} D_n}
|a_m||x^*(x_m)| &\leq
\sum_{m \in M_4 \setminus \bigcup_{ n \in I} D_n} |a_m| +
\sum_{m \in M_3 \setminus M_4} |x^*(x_m)| \notag \\
&< \epsilon + 
\sum_{ m \in M_3 \setminus \bigcup_{n=1}^{q}
       supp {\xi}_{i_n}^{P} }
 |x^*(x_m)| , \, \text{ by } (\ref{E:9}) \notag \\
&\text{ and since } \,
P \in [M_4], \notag \\
&< \epsilon + \epsilon = 2 \epsilon , \, \text{ by } (\ref{E:7}).
\end{align}
We also observe that
\begin{equation} \label{E:12}
\sum_{n \in I} |x^*(x_{l_{i_n}})| < \epsilon , \, \text{ as } \,
l_{i_n} \in 
M_3 \setminus \bigcup_{n=1}^{q}
       supp {\xi}_{i_n}^{P}  , \, 
\text{ for all } \, n \in I .
\end{equation}
Hence,
\begin{align}
x^* \biggl ( \sum_{m \leq m_0} a_m x_m \biggr ) &\geq
\sum_{ m \in \bigcup_{ n \in I} D_n} a_m x^*(x_m) -
\sum_{m \in M_3 \setminus \bigcup_{ n \in I} D_n}
|a_m||x^*(x_m)| \notag \\
&> \sum_{n \in I} \sum_{m \in D_n} a_m x^*(x_m) -
2 \epsilon , \, \text{ by } (\ref{E:11}), \notag \\
&\geq \sum_{n \in I} \sum_{m \in D_n \setminus \{l_{i_n}\}} 
a_m x^*(x_m) -
\sum_{n \in I} |x^*(x_{l_{i_n}})| - 2 \epsilon \notag \\
&\geq \sum_{n \in I} \sum_{m \in D_n \setminus \{l_{i_n}\}}
a_m \frac{1}{2} x^*({\xi}_{i_n}^{P} \cdot s)
- 3 \epsilon, \, \text{ by } \, (\ref{E:12}) \notag \\
&\text{ and since } \, D_n \setminus \{l_{i_n}\} \subset E_n, \notag \\
&\geq \frac{\delta}{4} 
\sum_{n \in I} x^*({\xi}_{i_n}^{P} \cdot s) - 3 \epsilon, \,
\text{ by } \, (\ref{E:10}) \notag \\
&\geq \frac{\delta}{4} ( \lambda - 2 \epsilon) - 3 \epsilon , \, 
\text{ by } \, (\ref{E:8}). \notag
\end{align}
It follows now that
\[ 1 \geq \biggl \| \sum_{m \leq m_0} a_m x_m \biggr \|
> \frac{\delta}{4} ( \lambda - 2 \epsilon) - 3 \epsilon ,\]
and thus \(\lambda < \frac{4}{\delta} (4 \epsilon +1)\)
contradicting the choice of \(\lambda\).

Hence, our claim holds and 
\(({\xi}_{i}^{P} \cdot s)\) is equivalent to the unit vector
basis of \(c_0\). Moreover, the equivalence constant \(C\),
depends only on \(K\), \(\delta\) and \(\epsilon\).
It is now easily seen that the set
\[ \{ P \in [M_5] : \, ({\xi}_{i}^{P} \cdot s) \, 
\text{ is }\, C 
\text{-equivalent to the } \, c_0 \text{-basis } \},
\]
is closed in \([M_5]\) and therefore Ramsey. Our previous argument 
yields \(N \in [M_5]\) so that
\(({\xi}_{i}^{Q} \cdot s)\) is equivalent to the unit vector
basis of \(c_0\), for all \(Q \in [N]\).
The proof of Theorem \ref{T:2} is now complete.
\end{proof}
Theorem \ref{T:1} follows immediately from Theorem \ref{T:2}.

\section{Non-trivial weak Cauchy sequences}
The last section is devoted to the relation between
Theorems \ref{T:1} and \ref{hr}. We first observe
the following immediate consequence of Theorem \eqref{hr}
\begin{Cor} \label{rc1}
The following are equivalent for a non-trivial weak Cauchy sequence
\((x_n)\) in a Banach space:

1. There exists a subsequence of \((x_n)\) which is \((s.s.)\).

2. There exists a subsequence \((x_{m_n})\) of \((x_n)\) such that
every convex block subsequence of \((x_{m_n})\) is semi-boundedly
complete.
\end{Cor}
\begin{proof}
The fact that 1. implies 2. is immediate since \((s.s.)\)
sequences are easily seen to be
semi-boundedly complete, and every convex block
subsequence of an \((s.s.)\) sequence is also \((s.s.)\)
\cite{ro}.

Suppose now that 2. holds. If no subsequence of \((x_n)\)
were \((s.s.)\), then Theorem \eqref{hr} yields a convex
block subsequence of \((x_{m_n})\) equivalent to the summing
basis. But the summing basis is not semi-boundedly complete.
This contradiction shows that 1. must hold.
\end{proof}
Corollary \ref{rc1} yields the following equivalent formulation
of Rosenthal's theorem:
\begin{Cor} \label{rc2}
For a non-trivial weak Cauchy sequence
\((x_n)\) in a Banach space one of the following statements
holds exclusively:

1. There exists a subsequence \((x_{m_n})\) of \((x_n)\)
such that every convex block subsequence of \((x_{m_n})\)
is semi-boundedly complete.

2. Every subsequence of \((x_n)\) admits a convex block 
subsequence equivalent to the summing basis.
\end{Cor}
Evidently, Corollary \ref{rc2} makes even more transparent
the analogy between Rosenthal's result and Theorem \ref{T:1}.
We also note here that, as it is shown in \cite{ro}, 
\((x_n)\) is an \((s.s.)\) sequence, if and only if
every proper subsequence of
its difference sequence is semi-boundedly complete.
(In the terminology of \cite{ro}, the difference sequence \((e_i)\)
of an \((s.s.)\) sequence is \((c.c.)\). That is, if
\(\sup_n \| \sum_{i=1}^n a_i e_i \|\) is finite, then
the scalar sequence \((a_n)\) converges.)

We next give a quantitative version of Theorem \ref{hr}.
\begin{Cor} \label{rc3}
Let \(t = (x_n)\) be a non-trivial weak Cauchy sequence having no
subsequence which is \((s.s.)\). Then for every
\(N \in [\mathbb{N}]\) there exist \(M \in [N]\),
a countable ordinal \(\xi\) and a constant \(C > 0\), so that
\(({\xi}_n^L \cdot t)\) is \(C\)-equivalent to the summing basis
for every \(L \in [M]\).
\end{Cor}
The proof of this corollary requires the following lemma
\begin{Lem} \label{rcl}
Let \(\xi\) be a countable ordinal, \(P \in [\mathbb{N}]\). Also let
\(t = (x_n)\) be a sequence in a Banach space which is 
\((\xi,L)\) convergent for every \(L \in [P]\).
Given \((\epsilon_n)\), a sequence of positive
scalars, there exists \(M \in [P]\) such that 
\(\|{\xi}_n^L \cdot t\| < \epsilon_n\), for every \(L \in [M]\)
and all \(n \in \mathbb{N}\).
\end{Lem}
\begin{proof}
Let \(Q \in [P]\). Our assumptions allow us
to choose \((k_n)\), an infinite subset
of \([\mathbb{N}]\), so that
\(\|{\xi}_{k_n}^Q \cdot t\| < \epsilon_n\), for all
\(n \in \mathbb{N}\). Stability now yields 
\(L \in [Q]\) such that
\(\|{\xi}_{n}^L \cdot t\| < \epsilon_n\), 
for all \(n \in \mathbb{N}\).
The assertion of the lemma follows from this since
\(\{ L \in [P] : \, \|{\xi}_{n}^L \cdot t\| < \epsilon_n , \,
\text{ for all } n \in \mathbb{N}\}\)
is a closed subset of \([P]\) and therefore Ramsey.
\end{proof} 
\begin{proof}[Proof of Corollary \ref{rc3}]
Theorem \eqref{hr} yields
\(u = (u_n)\) a convex block subsequence of 
\((x_n)\) equivalent to the summing basis.   
We set \(v= (x_n - u_n)\) which is clearly a weakly null sequence.
By the results of \cite{aa} there exists a countable ordinal
\(\xi\) such that no subsequence of \(v\) is an
\(\ell_1^{\xi}\) spreading model.
It follows now, by Theorem \ref{T:3}, that there
exists \(P \in [N]\) such that 
\(\lim_n \| {\xi}_n^L \cdot v \| = 0\), for every
\(L \in [P]\). 
We next choose \((\epsilon_n)\), a sequence of positive
scalars such that
\(\sum_n \epsilon_n < 1\).
As a consequence of Lemma \ref{rc3} we obtain 
\(M \in [P]\) such that
\(\|{\xi}_{n}^L \cdot v\| < \epsilon_n\), 
for every \(L \in [M]\)
and all \(n \in \mathbb{N}\). 
A standard perturbation result now yields \(D > 0\) such that   
\(({\xi}_{n}^L \cdot t)\)
is \(D\)-equivalent to \(({\xi}_{n}^L \cdot u)\),
for all \(L \in [M]\).
Since the summing basis is uniformly equivalent to all of its
convex block subsequences, we obtain \(C > 0\) such that 
\(({\xi}_n^L \cdot t)\) is \(C\)-equivalent to the summing
basis, for all \(L \in [M]\). This completes the proof.
\end{proof}  
We observe the similarity between the statements of
Theorem \ref{T:2} and Corollary \ref{rc3}.
However, the set of ordinals \(\xi\) satisfying the conclusion
of Theorem \ref{T:2} is a bounded segment of \([0, \omega_1)\),
in contrast with the corresponding set in Corollary \ref{rc3}
which is of course unbounded.
 
We continue our discussion about 
the relation between Theorems \ref{T:1} and \ref{hr}.
Recall that Rosenthal's \(c_0\)-theorem
states that every non-trivial weak Cauchy 
sequence \((x_j)\), either has
an \((s.s.)\) subsequence, or else there exists
a convex block subsequence \((s_j)\)
of \((x_j)\) equivalent to the summing basis.
In the later case, setting 
\(v_j = x_j - s_j\), we observe that
\((v_j)\) is weakly null and that
\((x_j - v_j) \) is equivalent to the
summing basis. Therefore, passing to a convex block 
subsequence in Rosenthal's theorem acts as a filtration
to remove the ``noise'' coming from an arbitrary
weakly null sequence. In our case the reasoning for
passing to a convex block subsequence is different:
It exhausts the local \(\ell_1\) structure of the sequence.

In spite these differences it seems that
there are similarities in the statements for weakly null
and non-trivial weak Cauchy sequences. Our final
corollary which is the analog to Elton's dichotomy illustrates this
\begin{Cor}
For a non-trivial weak Cauchy sequence \((x_n)\) one of
the following statements holds exclusively:

1. Every subsequence of \((x_n)\) admits a subsequence
equivalent to the summing basis.

2. There exists a subsequence of \((x_n)\) which is
semi-boundedly complete.
\end{Cor}
\begin{proof}
Clearly the statements are mutually exclusive since the
summing basis is not semi-boundedly complete.

Suppose that 2. does not hold. It follows that no subsequence
of \((x_n)\) is \((s.s.)\). Theorem \ref{hr} now yields
that every subsequence of \((x_n)\) admits a convex block
subsequence equivalent to the summing basis.
To prove that 1. holds let \((x_{m_n})\) be a subsequence
of \((x_n)\). Without loss of generality, by passing to a
subsequence according to Proposition 2.2 of \cite{ro}, we can assume
that \((x_{m_n})\) dominates the summing basis.

Next choose \((u_n)\), a convex block subsequence of \((x_n)\)
equivalent to the summing basis. We set
\(y_n = x_{m_n} - u_n\), for all \(n \in \mathbb{N}\).
If \((y_n)\) is not semi-normalized, then 1. follows.
So assuming that \((y_n)\) is semi-normalized we  
claim that there exists a subsequence of \((y_n)\) equivalent
to the unit vector basis of \(c_0\). 
Indeed, if that were not the case, then by Elton's dichotomy,
Corollary \ref{T:33},
there would exist a semi-boundedly complete subsequence of \((y_n)\).
But since every subsequence of \((x_{m_n})\) dominates
the summing basis (and therefore every subsequence of \((u_n)\) as well),
we obtain that \((x_{m_n})\) has a 
semi-boundedly complete subsequence which of course contradicts
our assumption that 2. does not hold.
Hence, our claim must hold and it immediately yields
a subsequence of \((x_{m_n})\) equivalent to the summing
basis in view of the following elementary fact: Let \((f_n)\)
and \((g_n)\) be sequences in a Banach space with
\((f_n)\) equivalent to the summing basis and \((g_n)\)
equivalent to the \(c_0\) basis. Then there exists a subsequence of
\((f_n + g_n)\) equivalent to the summing basis. 
\end{proof}

\end{document}